\documentclass{amsart}%
\usepackage{amsfonts}
\usepackage{amsmath}
\usepackage{amssymb}
\usepackage{graphicx}%
\providecommand{\U}[1]{\protect\rule{.1in}{.1in}}
\newtheorem{theorem}{Theorem}
\theoremstyle{plain}

\newtheorem{corollary}{Corollary}

\newtheorem{definition}{Definition}
\newtheorem{example}{Example}

\newtheorem{proposition}{Proposition}
\newtheorem{remark}{Remark}

\newtheorem{summary}{Summary}
\numberwithin{equation}{section}
\begin{document}
\title[Symbolic Number Theory]{Symbolic Arithmetic and Integer Factorization}
\author{Samuel J. Lomonaco}
\address{University of Maryland Baltimore County (UMBC)\\
Baltimore, MD \ 21250\\
and Princeton University, Princeton, NJ \ 08540}
\email{Lomonaco@umbc.edu}
\urladdr{http://www.csee.umbc.edu/\symbol{126}lomonaco}
\thanks{The author would like to thank Gamal Abdali and Sumeet Bagde for implementing
the BF algorithm respectively in LISP and in Mathematica. Thanks are also due
to John Conway, Alan Sherman, and Katie Smith whose supportive interest in
this paper encouraged the author to write up these results. \ Thanks are also
due to the Princeton University Mathematics Department for providing, during
my sabbatical visit, a congenial, encouraging research atmosphere in which to
complete this paper.}
\date{March 10, 2013}

\begin{abstract}
In this paper, we create a systematic and automatic procedure for transforming
the integer factorization problem into the problem of solving a system of
Boolean equations. \ Surprisingly, the resulting system of Boolean equations
takes on a "life of its own" and becomes a new type of integer, which we call
a \textbf{generic integer}. \ \newline We then proceed to use the newly found
algebraic structure of the\textbf{ ring of generic integers} to create two new
integer factoring algorithms, called respectively the \textbf{Boolean
factoring (BF) algorithm}, and the\textbf{ multiplicative Boolean factoring
(MBF) algorithm}. \ Although these two algorithms are not competitive with
current classical integer factoring algorithms, it is hoped that they will
become stepping stones to creating much faster and more competitive
algorithms, and perhaps be precursors of a new quantum algorithm for integer factoring.

\end{abstract}
\maketitle
\tableofcontents

\section{Introduction}

This paper is the result of a research program seeking to gain a better
qualitative and quantitative understanding of the computational complexity of
integer factorization. \ By the phrase \textquotedblleft computational
complexity\textquotedblright\ is meant the \textbf{Boolean complexity}, i.e.,
the minimum number of fundamental Boolean operations required to factor an
integer, as a function of integer size. \ The strategy chosen for
accomplishing this research objective was to develop a systematic and
automatic procedure for the decomposition of arithmetic operations into
Boolean operations. \ 

\bigskip

The results of this endeavor produced unanticipated results. \ As expected,
what resulted was a conversion of the problem of integer factoring into the
problem of solving a system of Boolean equations. \ But these resulting
systems of Boolean equations unexpectedly took on a "life of their own," and
became a new type of integer in their own right. We call this new type of
integer a \textbf{generic integer}, and the corresponding algebraic object
$\mathbb{G}\left\langle \underline{\underline{x}}\right\rangle $, the
\textbf{ring of generic integers}. \ Such a ring is a fascinating mix of
characteristic $0$ and characteristic $2$ algebraic structure.

\bigskip

Yet another surprise was that this symbolic approach to integer factoring
naturally led to the larger context of the ring of \textbf{dyadic integers}
$\mathbb{Z}_{(2)}$, and then on to the corresponding ring of \textbf{generic
dyadic integers} $\mathbb{G}_{(2)}\left\langle \underline{\underline{x}%
}\right\rangle $. \ The advantage of looking at integer factorization within
the larger context of the dyadics is that every odd (generic) integer dyadic
has a (generic) dyadic integer inverse.\ \ This naturally led to the creation
of the (\textbf{generic}) \textbf{lopsided division algorithm} for computing
inverses of odd integers in the dyadic integers $Z_{(2)}$ (in the generic
dyadic integers $\mathbb{G}_{(2)}\left\langle \underline{\underline{x}%
}\right\rangle $).

\bigskip

This in turn led to the creation of the \textbf{Boolean factoring (BF)
algorithm}, which systematically and automatically translates the problem of
factoring an integer $N$ into the problem of solving a system of Boolean
equations. \ These Boolean equations were obtained by using generic lopsided
division to divide the integer $N$ by carefully chosen odd generic integer
$x$. \ The sought system of Boolean equations is simply the generic remainder
resulting from this division.

\bigskip

Next it is noted that the generic inverse $x^{-1}$ of $x$ need only be
computed once, and then used over and over again. \ This immediately leads to
the creation of a second factoring algorithm, the \textbf{multiplicative
Boolean factoring (MBF) algorithm}, which simply computes the generic product%
\[
N\cdot x^{-1}\text{ ,}%
\]
to produce the system of Boolean equations to be solved.

\bigskip

A general framework (i.e., scarce satisfiability) for solving the system of
Boolean equations produced by the BF and MBF algorithms is then discussed.
\ This framework was later used by Gamal Abdali to create a LISP
implementation of the BF algorithm. \ Sumeet Bagde then extended these methods
by using binary decision diagrams (BDDs)\cite{Bryant2} to create a Mathematica
program that also implemented the BF algorithm.

\bigskip

Both the LISP and Mathematica programs were used to factor a large number of
integers. \ The runtime statistics indicated that the BF algorithm runs in
exponential time, and hence, is not competitive with the best classical
factoring algorithms. \ For an algebraic basis as to why this is the case, we
refer the reader to the topdown overview given toward the end of this paper.

\bigskip

Open questions and future possible research directions are discussed in the
conclusion of this paper. \ Connections with the satisfiability problem are
also discussed.

\bigskip

\section{Lopsided division}

Let $\mathbb{Z}_{(2)}$ denote the ring of dyadic integers, and let
$\mathbb{Z}$ denote its subring of all rational integers (i.e., its subring of
all standard integers.) \bigskip\ 

Given below are examples of the dyadic expansion of some rational integers.
Please note that the dyadic expansion of a non-negative rational integer is
the conventional radix 2 expansion. The dyadic expansion of a negative
rational integer is an \textquotedblleft infinite 2's
complement\textquotedblright\ of the radix 2 expansion of its absolute value.
\[%
\begin{array}
[c]{rcc}%
5 & = & \ldots\;00101\\
&  & \\
4 & = & \ldots\;00100\\
&  & \\
3 & = & \ldots\;00011\\
&  & \\
2 & = & \ldots\;00010\\
&  & \\
1 & = & \ldots\;00001\\
&  & \\
0 & = & \ldots\;00000\\
&  & \\
-1 & = & \ldots\;11111\\
&  & \\
-2 & = & \ldots\;11110\\
&  & \\
-3 & = & \ldots\;11101\\
&  & \\
-4 & = & \ldots\;11100\\
&  & \\
-5 & = & \ldots\;11011
\end{array}
\]

\bigskip\ 

Every odd dyadic integer is a unit in the ring of dyadic integers
$\mathbb{Z}_{(2)}$ , i.e., given any odd dyadic integer $a$ , there exists a
unique dyadic integer $a^{-1}$ such that
\[
a\cdot a^{-1}=1
\]

Thus,

\begin{proposition}
Let $b$ be an odd dyadic integer, and let $a$ be an arbitrary dyadic integer.
Then $a$ \textbf{divided by} $b$ , written
\[
a/b
\]
is also a dyadic integer.
\end{proposition}

\bigskip\ 

\begin{corollary}
Let $b$ be an odd rational integer, and let $a$ be an arbitrary rational
integer. Then $a$ \textbf{divided by} $b$ , written
\[
a/b
\]
is a well defined dyadic integer. Moreover, the dyadic integer $a/b$ is a
rational integer if and only if $b$ is an exact divisor of $a$ in the ring of
rational integers $\mathbb{Z}$.
\end{corollary}

\bigskip\ 

\begin{definition}
Let
\[
\ldots\;,\;a^{(2)}\;,\;a^{(1)}\;,\;a^{(0)}
\]
be a sequence of dyadic integers, and let
\[
a_{j}^{(i)}
\]
denote the $j$-th bit of the dyadic expansion of $a^{(i)}$. Then the sequence
\[
\ldots\;,\;a^{(2)}\;,\;a^{(1)}\;,\;a^{(0)}
\]
is said to be \textbf{convergent} provided for each $j\geq0$ there exists a
non-negative integer $n=n(j)$ such that
\[
a_{j}^{(i)}=a_{j}^{(n(j))}\quad\text{for}\quad i\geq n(j).
\]
Otherwise, the sequence is said to be \textbf{divergent}. If the above
sequence is convergent, its \textbf{limit}, written
\[
\lim\limits_{i\rightarrow\infty}a^{(i)}
\]
is said to \textbf{exist}, and is defined as the dyadic integer with dyadic
expansion given by
\[
\left(  \lim\limits_{i\rightarrow\infty}a^{(i)}\right)  _{j}=a_{j}^{(n(j))}.
\]

\end{definition}

\bigskip\ 

\begin{remark}
The above limit is the standard limit in the valuation topology.
\end{remark}

\bigskip\ 

\begin{definition}
Let $a$ be a dyadic integer with dyadic expansion
\[
\ldots\;,\;a_{2}\;,\;a_{1}\;,\;a_{0}\quad.
\]
Then let
\[
\mathcal{S}:\mathbb{Z}_{(2)}\longrightarrow\mathbb{Z}_{(2)}%
\]
denote a \textbf{left shift} by 1 bit, i.e.,
\[
\mathcal{S}a=\ldots\;,\;a_{2}\;,\;a_{1}\;,\;a_{0}\;,\;0
\]

\end{definition}

\bigskip\ 

\begin{remark}
Hence, $\mathcal{S}a$ is the same as $a$ multiplied by the dyadic
\[
2=\ldots\;0\;0\;0\;1\;0
\]

\end{remark}

\bigskip\ 

\begin{definition}
Let $a$, $b$, and $c$ be dyadic integers with dyadic expansions
\[
\left\{
\begin{array}
[c]{c}%
\ldots\;,\;a_{2}\;,\;a_{1}\;,\;a_{0}\\
\\
\ldots\;,\;b_{2}\;,\;b_{1}\;,\;b_{0}\\
\\
\ldots\;,\;c_{2}\;,\;c_{1}\;,\;c_{0}%
\end{array}
\right.
\]
respectively. Define the \textbf{first bitwise symmetric function} of $a$,
$b$, $c$, written
\[
Bitwise\_\sigma_{1}(a,b,c)\;,
\]
as the dyadic with $i$-th bit given by
\[
\sigma_{1}(a_{i},b_{i},c_{i})=a_{i}\dotplus b_{i}\dotplus c_{i}%
\]
where \textquotedblleft$\dotplus$\textquotedblright\ denotes the exclusive
\textquotedblleft or\textquotedblright\ binary operation. Define the
\textbf{second bitwise symmetric function} of $a$, $b$, $c$, written
\[
Bitwise\_\sigma_{2}(a,b,c)\;,
\]
as the dyadic with $i$-th bit given by
\[
\sigma_{2}(a_{i},b_{i},c_{i})=\left(  a_{i}\diamond b_{i}\right)
\dotplus\left(  b_{i}\diamond c_{i}\right)  \dotplus\left(  c_{i}\diamond
a_{i}\right)  \;,
\]
where \textquotedblleft$\dotplus$\textquotedblright\ again denotes exclusive
\textquotedblleft or\textquotedblright\ and \textquotedblleft$\diamond
$\textquotedblright\ denotes logical \textquotedblleft and\textquotedblright\ .
\end{definition}

\bigskip\ 

We are now ready to define an algorithm called \textbf{lopsided division}.

\bigskip\ 

\begin{theorem}
Let $b$ be an odd dyadic integer, and let $a$ be an arbitrary dyadic integer,
with dyadic expansions
\[
\left\{
\begin{array}
[c]{c}%
\ldots\;,\;b_{2}\;,\;b_{1}\;,\;b_{0}\\
\\
\ldots\;,\;a_{2}\;,\;a_{1}\;,\;a_{0}%
\end{array}
\right.  ,
\]
respectively (where $b_{0}=1$.) Let
\[
\left\{
\begin{array}
[c]{lll}%
c^{(0)} & = & a\\
&  & \\
borrows^{(0)} & = & 0
\end{array}
\right.
\]
and
\[
\left\{
\begin{array}
[c]{lll}%
c^{(i+1)} & = & Bitwise\_\sigma_{1}\left(  \;c^{(i)}\;,borrows^{(i)}%
\;,\;c_{i}^{(i)}\diamond\mathcal{S}^{i}b\;\right) \\
&  & \\
borrows^{(i+1)} & = & \mathcal{S}\cdot Bitwise\_\sigma_{2}\left(
\;c^{(i)\ast}\;,borrows^{(i)}\;,\;c_{i}^{(i)}\diamond\mathcal{S}%
^{i}b\;\right)
\end{array}
\right.
\]
where $c^{(i)\ast}$ denotes the \textbf{complement} of $c^{(i)}$ , i.e.,
\[
c^{(i)\ast}=1\dotplus c^{(i)}\qquad\qquad\qquad\text{for }j\geq0
\]
where `$\dotplus$' denotes exclusive \textquotedblleft or\textquotedblright,
and where $c_{i}^{(i)}\diamond\mathcal{S}^{i}b$ denotes the bitwise logical
\textquotedblleft and\textquotedblright, i.e.,
\[
\left(  c_{i}^{(i)}\diamond\mathcal{S}^{i}b\right)  _{j}=c_{i}^{(i)}%
\diamond\left(  \mathcal{S}^{i}b\right)  _{j}%
\]
Then the sequences $c^{(i)}$ and $borrows^{(i)}$ are convergent, converging
respectively to:
\[
\left\{
\begin{array}
[c]{lll}%
\lim\limits_{i\rightarrow\infty}c^{(i)} & = & a/b\\
&  & \\
\lim\limits_{i\rightarrow\infty}borrows^{(i)} & = & 0
\end{array}
\right.
\]
This algorithm for computing $a/b$ is called \textbf{lopsided division}.
\end{theorem}

\bigskip\ 

\begin{remark}
Please note that this is an algorithm in the sense that
\[
\left(  a/b\right)  _{j}=c_{j}^{(i)}%
\]

\end{remark}

\bigskip\ 

\begin{definition}
Let $a$ be a positive integer. The \textbf{length} of $a$, written
\[
lgth(a)\qquad,
\]
is defined as
\[
lgth(a)=j+1
\]
where $j$ is the largest non-negative integer such that $a_{j}=1$.
\end{definition}

\bigskip\ 

\begin{corollary}
Let $b$ be a positive odd rational integer, and let $a$ be an arbitrary
positive rational integer. Let
\[
\Gamma=1+lgth(a)-lgth(b)
\]
Then $b$ is an exact divisor of $a$ if and only if
\[
c^{(\Gamma)}=borrows^{(\Gamma)}\qquad.
\]
If $b$ is an exact divisor of $a$, then the radix 2 expansion of $a/b$ is
\[
a/b=c_{\Gamma-1}^{(\Gamma-1)}\quad c_{\Gamma-2}^{(\Gamma-2)}\quad\cdots\quad
c_{1}^{(1)}\quad c_{0}^{(0)}\quad.
\]

\end{corollary}

\bigskip\ 

\begin{example}
Lopsided division of $a=209$ by $b=19$ . The radix 2 representations of $209$
and $19$ are respectively:
\[
\left\{
\begin{array}
[c]{lll}%
a & = & 11010001\\
&  & \\
b & = & 10011
\end{array}
\right.
\]
The lopsided division of $209$ by $19$ is given in the tableau below:
\[%
\begin{tabular}
[c]{ccccccccccc}
&  &  &  &  & $c_{3}^{(3)}$ & $c_{2}^{(2)}$ & $c_{1}^{(1)}$ & $c_{0}^{(0)}$ &
& \\
&  &  &  &  & $1$ & $0$ & $1$ & $1$ & $a/b$ & \\\cline{2-9}\cline{9-9}%
$a=$ & $\ 1\ $ & $\ 1\ $ & $\ 0\ $ & $\ 1\ $ & $\ 0\ $ & $\ 0\ $ & $\ 0\ $ &
$\ 1\ $ & \multicolumn{1}{|c}{$10011=b$} & \\
&  &  &  & $\ 1\ $ & $\ 0\ $ & $\ 0\ $ & $\ 1\ $ & $\ 1\ $ & $c_{0}%
^{(0)}\diamond\mathcal{S}^{0}b$ & \\\cline{5-9}
&  &  &  & $\ 0\ $ & $\ 0\ $ & $\ 0\ $ & $\ 1\ $ & $\ 0\ $ & $c^{(1)}$ & \\
&  &  & $\ 0\ $ &  &  & $\ 1\ $ &  &  & $borrows^{(1)}$ & \\
&  &  & $\ 1\ $ & $\ 0\ $ & $\ 0\ $ & $\ 1\ $ & $\ 1\ $ &  & $c_{1}%
^{(1)}\diamond\mathcal{S}^{1}b$ & \\\cline{4-9}
&  &  & $\ 1\ $ & $\ 0\ $ & $\ 0\ $ & $\ 0\ $ & $\ 0\ $ &  & $c^{(2)}$ & \\
&  & $\ 1\ $ &  &  & $\ 1\ $ &  &  &  & $borrows^{(2)}$ & \\
&  & $\ 0\ $ & $\ 0\ $ & $\ 0\ $ & $\ 0\ $ & $\ 0\ $ &  &  & $c_{2}%
^{(2)}\diamond\mathcal{S}^{2}b$ & \\\cline{3-7}
&  & $\ 0\ $ & $\ 1\ $ & $\ 0\ $ & $\ 1\ $ & $\ 0\ $ &  &  & $c^{(3)}$ & \\
& $\ 0\ $ &  &  & $\ 1\ $ &  &  &  &  & $borrows^{(3)}$ & \\
& $\ 1\ $ & $\ 0\ $ & $\ 0\ $ & $\ 1\ $ & $\ 1\ $ &  &  &  & $c_{3}%
^{(3)}\diamond\mathcal{S}^{3}b$ & \\\cline{2-6}
& $\ 0\ $ & $\ 0\ $ & $\ 1\ $ & $\ 0\ $ & $\ 0\ $ &  &  &  & $c^{(4)}$ & \\
&  &  & $\ 1\ $ &  &  &  &  &  & $borrows^{(4)}$ &
\end{tabular}
\ \ \ \
\]
Please note that
\[%
\begin{array}
[c]{lll}%
1+lgth(209)-lgth(19) & = & 1+8-5=4\\
& \text{and} & \\
c^{(4)} & = & borrows^{(4)}%
\end{array}
\]
Hence, $19$ is an exact divisor of $209$, and
\[
209/19=11\;(base\;10)=10011\;(base\;2).
\]

\end{example}

\bigskip\ 

\begin{example}
Lopsided division of $209$ by $21$.
\[%
\begin{tabular}
[c]{ccccccccccc}
&  &  &  &  & $c_{3}^{(3)}$ & $c_{2}^{(2)}$ & $c_{1}^{(1)}$ & $c_{0}^{(0)}$ &
& \\
&  &  &  &  & $1$ & $0$ & $1$ & $1$ & $a/b$ & \\\cline{2-9}\cline{9-9}%
$a=$ & $\ 1\ $ & $\ 1\ $ & $\ 0\ $ & $\ 1\ $ & $\ 0\ $ & $\ 0\ $ & $\ 0\ $ &
$\ 1\ $ & \multicolumn{1}{|c}{$10101=b$} & \\
&  &  &  & $\ 1\ $ & $\ 0\ $ & $\ 1\ $ & $\ 0\ $ & $\ 1\ $ & $c_{0}%
^{(0)}\diamond\mathcal{S}^{0}b$ & \\\cline{5-9}
&  &  &  & $\ 0\ $ & $\ 0\ $ & $\ 1\ $ & $\ 0\ $ & $\ 0\ $ & $c^{(1)}$ & \\
&  &  & $\ 0\ $ &  & $\ 1\ $ &  &  &  & $borrows^{(1)}$ & \\
&  &  & $\ 0\ $ & $\ 0\ $ & $\ 0\ $ & $\ 0\ $ & $\ 0\ $ &  & $c_{1}%
^{(1)}\diamond\mathcal{S}^{1}b$ & \\\cline{4-9}
&  &  & $\ 0\ $ & $\ 0\ $ & $\ 1\ $ & $\ 1\ $ & $\ 0\ $ &  & $c^{(2)}$ & \\
&  & $\ 0\ $ &  &  & $\ 1\ $ &  &  &  & $borrows^{(2)}$ & \\
&  & $\ 1\ $ & $\ 0\ $ & $\ 1\ $ & $\ 0\ $ & $\ 1\ $ &  &  & $c_{2}%
^{(2)}\diamond\mathcal{S}^{2}b$ & \\\cline{3-7}
&  & $\ 0\ $ & $\ 0\ $ & $\ 0\ $ & $\ 1\ $ & $\ 0\ $ &  &  & $c^{(3)}$ & \\
& $\ 0\ $ &  & $\ 1\ $ &  &  &  &  &  & $borrows^{(3)}$ & \\
& $\ 1\ $ & $\ 0\ $ & $\ 1\ $ & $\ 0\ $ & $\ 1\ $ &  &  &  & $c_{3}%
^{(3)}\diamond\mathcal{S}^{3}b$ & \\\cline{2-6}
& $\ 0\ $ & $\ 0\ $ & $\ 0\ $ & $\ 0\ $ & $\ 0\ $ &  &  &  & $c^{(4)}$ & \\
&  &  & $\ 1\ $ &  &  &  &  &  & $borrows^{(4)}$ &
\end{tabular}
\ \ \ \ \
\]
Please note that
\[%
\begin{array}
[c]{lll}%
1+lgth(209)-lgth(21) & = & 1+8-5=4\\
& \text{and} & \\
c^{(4)} & \neq & borrows^{(4)}%
\end{array}
\]
Hence, $a/b=209/21$ is not a rational integer, and $b=21$ is not an exact
divisor of $a=209$. So
\[%
\begin{array}
[c]{ccc}%
1101\;(base\;2) & \neq & \left\lfloor 209/21\right\rfloor
\end{array}
\]

\end{example}

\bigskip\ 

\begin{example}
Lopsided division of $209$ by $17$.
\[%
\begin{tabular}
[c]{ccccccccccc}
&  &  &  &  & $c_{3}^{(3)}$ & $c_{2}^{(2)}$ & $c_{1}^{(1)}$ & $c_{0}^{(0)}$ &
& \\
&  &  &  &  & $0$ & $0$ & $0$ & $1$ & $a/b$ & \\\cline{2-9}\cline{9-9}%
$a=$ & $\ 1\ $ & $\ 1\ $ & $\ 0\ $ & $\ 1\ $ & $\ 0\ $ & $\ 0\ $ & $\ 0\ $ &
$\ 1\ $ & \multicolumn{1}{|c}{$10001=b$} & \\
&  &  &  & $\ 1\ $ & $\ 0\ $ & $\ 0\ $ & $\ 0\ $ & $\ 1\ $ & $c_{0}%
^{(0)}\diamond\mathcal{S}^{0}b$ & \\\cline{5-9}
&  &  &  & $\ 0\ $ & $\ 0\ $ & $\ 0\ $ & $\ 0\ $ & $\ 0\ $ & $c^{(1)}$ & \\
&  &  & $\ 0\ $ &  &  &  &  &  & $borrows^{(1)}$ & \\
&  &  & $\ 0\ $ & $\ 0\ $ & $\ 0\ $ & $\ 0\ $ & $\ 0\ $ &  & $c_{1}%
^{(1)}\diamond\mathcal{S}^{1}b$ & \\\cline{4-9}
&  &  & $\ 0\ $ & $\ 0\ $ & $\ 0\ $ & $\ 0\ $ & $\ 0\ $ &  & $c^{(2)}$ & \\
&  & $\ 0\ $ &  &  &  &  &  &  & $borrows^{(2)}$ & \\
&  & $\ 0\ $ & $\ 0\ $ & $\ 0\ $ & $\ 0\ $ & $\ 0\ $ &  &  & $c_{2}%
^{(2)}\diamond\mathcal{S}^{2}b$ & \\\cline{3-7}
&  & $\ 1\ $ & $\ 0\ $ & $\ 0\ $ & $\ 0\ $ & $\ 0\ $ &  &  & $c^{(3)}$ & \\
& $\ 0\ $ &  &  &  &  &  &  &  & $borrows^{(3)}$ & \\
& $\ 0\ $ & $\ 0\ $ & $\ 0\ $ & $\ 0\ $ & $\ 0\ $ &  &  &  & $c_{3}%
^{(3)}\diamond\mathcal{S}^{3}b$ & \\\cline{2-6}
& $\ 1\ $ & $\ 1\ $ & $\ 0\ $ & $\ 0\ $ & $\ 0\ $ &  &  &  & $c^{(4)}$ & \\
$\ 0\ $ &  &  &  &  &  &  &  &  & $borrows^{(4)}$ &
\end{tabular}
\ \ \ \
\]
Please note that
\[%
\begin{array}
[c]{lll}%
1+lgth(209)-lgth(17) & = & 1+8-5=4\\
& \text{and} & \\
c^{(4)} & \neq & borrows^{(4)}%
\end{array}
\]
Hence, $17$ is not an exact divisor of $209$.
\end{example}

\bigskip\ 

\begin{example}
Lopsided division of $513$ by $27$.
\[%
\begin{tabular}
[c]{ccccccccccccc}
&  &  &  &  & $c_{5}^{(5)}$ & $c_{4}^{(4)}$ & $c_{3}^{(3)}$ & $c_{2}^{(2)}$ &
$c_{1}^{(1)}$ & $c_{0}^{(0)}$ &  & \\
&  &  &  &  & $0$ & $1$ & $0$ & $0$ & $1$ & $1$ & $a/b$ & \\\cline{2-11}%
\cline{11-11}%
$a=$ & $\ 1\ $ & $\ 0\ $ & $\ 0\ $ & $\ 0\ $ & $\ 0\ $ & $\ 0\ $ & $\ 0\ $ &
$\ 0\ $ & $\ 0\ $ & $\ 1\ $ & \multicolumn{1}{|c}{$11011=b$} & \\
&  &  &  &  &  & $\ 1\ $ & $\ 1\ $ & $\ 0\ $ & $\ 1\ $ & $\ 1\ $ &
$c_{0}^{(0)}\diamond\mathcal{S}^{0}b$ & \\\cline{7-11}
&  &  &  &  &  & $\ 1\ $ & $\ 1\ $ & $\ 0\ $ & $\ 1\ $ & $\ 0\ $ & $c^{(1)}$ &
\\
&  &  &  &  & $\ 1\ $ & $\ 1\ $ &  & $\ 1\ $ &  &  & $borrows^{(1)}$ & \\
&  &  &  &  & $\ 1\ $ & $\ 1\ $ & $\ 0\ $ & $\ 1\ $ & $\ 1\ $ &  &
$c_{1}^{(1)}\diamond\mathcal{S}^{1}b$ & \\\cline{6-11}
&  &  &  &  & $\ 0\ $ & $\ 1\ $ & $\ 1\ $ & $\ 0\ $ & $\ 0\ $ &  & $c^{(2)}$ &
\\
&  &  &  & $\ 1\ $ & $\ 1\ $ &  & $\ 1\ $ &  &  &  & $borrows^{(2)}$ & \\
&  &  &  & $\ 0\ $ & $\ 0\ $ & $\ 0\ $ & $\ 0\ $ & $\ 0\ $ &  &  &
$c_{2}^{(2)}\diamond\mathcal{S}^{2}b$ & \\\cline{5-9}
&  &  &  & $\ 1\ $ & $\ 1\ $ & $\ 1\ $ & $\ 0\ $ & $\ 0\ $ &  &  & $c^{(3)}$ &
\\
&  &  & $\ 1\ $ & $\ 1\ $ &  &  &  &  &  &  & $borrows^{(3)}$ & \\
&  &  & $\ 0\ $ & $\ 0\ $ & $\ 0\ $ & $\ 0\ $ & $\ 0\ $ &  &  &  &
$c_{3}^{(3)}\diamond\mathcal{S}^{3}b$ & \\\cline{4-8}
&  &  & $\ 1\ $ & $\ 0\ $ & $\ 1\ $ & $\ 1\ $ & $\ 0\ $ &  &  &  & $c^{(4)}$ &
\\
&  & $\ 1\ $ &  &  &  &  &  &  &  &  & $borrows^{(4)}$ & \\
&  & $\ 1\ $ & $\ 1\ $ & $\ 0\ $ & $\ 1\ $ & $\ 1\ $ &  &  &  &  & $c$%
$_{4}^{(4)}\diamond\mathcal{S}^{4}$$b$ & \\\cline{3-7}
&  & $\ 0\ $ & $\ 0\ $ & $\ 0\ $ & $\ 0\ $ & $\ 0\ $ &  &  &  &  & $c^{(5)}$ &
\\
& $\ 1\ $ &  &  &  &  &  &  &  &  &  & $borrows^{(5)}$ & \\
& $\ 0\ $ & $\ 0\ $ & $\ 0\ $ & $\ 0\ $ & $\ 0\ $ &  &  &  &  &  & $c$%
$_{5}^{(5)}\diamond\mathcal{S}^{5}$$b$ & \\\cline{2-6}
& $\ 0\ $ & $\ 0\ $ & $\ 0\ $ & $\ 0\ $ & $\ 0\ $ &  &  &  &  &  & $c$$^{(6)}$
& \\
$\ 0\ $ &  &  &  &  &  &  &  &  &  &  & $borrows$$^{(6)}$ &
\end{tabular}
\ \ \
\]

\end{example}

\bigskip\ 

Please note that
\[%
\begin{array}
[c]{lll}%
1+lgth(513)-lgth(27) & = & 1+10-5=6\\
& \text{and} & \\
c^{(6)} & = & borrows^{(6)}%
\end{array}
\]
Hence, $27$ is an exact divisor of $513$ . Moreover
\[
513/27=19=10011\;(base\;2)
\]

\bigskip\ 

\section{Generic arithmetic: ``Algebraic parallel processing''}

We now lift the algebraic operations on the ring of rational integers
$\mathbb{Z}$ and on the ring of dyadics integers $\mathbb{Z}_{(2)}$ to the
generic level by creating respectively the ring $\mathbb{G}\left\langle
\mathbf{x}\right\rangle $ of generic rational integers and the ring
$\mathbb{G}_{(2)}\left\langle \mathbf{x}\right\rangle $ of generic dyadic
integers. \ It will then be observed that each arithmetic operation in these
generic rings, $\mathbb{G}\left\langle \mathbf{x}\right\rangle $ and
$\mathbb{G}_{(2)}\left\langle \mathbf{x}\right\rangle $, is equivalent to
performing many simultaneous arithmetic operations in the corresponding
respective rings of rational integers $\mathbb{Z}$ and dyadic integers
$\mathbb{Z}_{(2)}$. \ This is what is meant by the phrase \textquotedblleft
algebraic parallel processing.\textquotedblright\ \ (Please refer to section 7
for a topdown overview.)

\bigskip

\begin{definition}
A \textbf{Boolean ring} $\mathbb{B}$ (with addition denoted by
\textquotedblleft$\dotplus$\textquotedblright\ and with multiplication denoted
by \textquotedblleft$\diamond$\textquotedblright\ ) is a ring with
multiplicative identity, denoted by \textquotedblleft$1$\textquotedblright,
such that each element $a$ of $\mathbb{B}$ is an idempotent, i.e., such
that$_{{}}$
\[
a^{2}=a.
\]
It follows that $\mathbb{B}$ is a commutative ring, and that each element $a$
of $\mathbb{B}$ is its own additive inverse, i.e.,
\[
a\dotplus a=0,
\]
where \textquotedblleft$0$\textquotedblright\ denotes the additive identity of
$\mathbb{B}$ . The additive operation will often be referred to as
\textbf{exclusive \textquotedblleft or\textquotedblright} and the
multiplicative operation \textquotedblleft$\diamond$\textquotedblright\ will
often be referred to as \textbf{logical \textquotedblleft
and\textquotedblright}. The \textbf{complement} of an element $a$ of
$\mathbb{B}$ , written $a^{\ast}$, is defined as
\[
a^{\ast}=1\dotplus a.
\]

\end{definition}

\bigskip\ 

\begin{definition}
Let $\nu$ be an arbitrary but fixed non-negative integer. \ Let
\[
x_{\nu-1},...,x_{2},x_{1},x_{0}%
\]
denote a finite sequence of $\nu$ distinct symbols, let $\mathbf{x}$ denote
the set of these symbols, and let
\[
\mathbb{B}\left\langle \mathbf{x}\right\rangle
\]
denote the \textbf{free Boolean ring} on the symbols in $\mathbf{x}$ . The
elements of $\mathbf{x}$ are called the \textbf{free basis elements} of
$\mathbb{B}<\mathbf{x}>$ and $\mathbf{x}$ is called the \textbf{free basis} of
$\mathbb{B}<\mathbf{x}>$.
\end{definition}

\bigskip\ 

\begin{remark}
Thus, $\mathbb{B}\left\langle \ \right\rangle $ denotes the free Boolean ring
on the empty free basis. Hence, $\mathbb{B}\left\langle \ \right\rangle $ may
be identified with the field of two elements $\mathbb{F}_{2}$.
\end{remark}

\bigskip\ 

A \textbf{generic dyadic integer} $e$ is a an infinite sequence
\[
e=\;\ldots,\;e_{2},\;e_{1},\;e_{0}%
\]
of elements $e_{i}$ of $\mathbb{B}\left\langle \mathbf{x}\right\rangle $ . If
there exists an integer $k$ such that
\[
e_{i}=e_{k}%
\]
for all $i\geq k$ , then $e$ is called a \textbf{generic rational integer}.

Let
\[
\mathbb{G}_{(2)}=\mathbb{G}_{(2)}\left\langle \mathbf{x}\right\rangle
=\mathbb{G}_{(2)}\left\langle x_{\nu-1},\;\ldots,\;x_{2},\;x_{1}%
,\;x_{0}\right\rangle
\]
and
\[
\mathbb{G}=\mathbb{G}\left\langle \mathbf{x}\right\rangle =\mathbb{G}%
\left\langle x_{\nu-1},\;\ldots\;,\;x_{2},\;x_{1},\;x_{0}\right\rangle
\]
denote respectively the set of all generic dyadic integers and the set of all
generic rational integers.

The generic integers
\[
\left\{
\begin{array}
[c]{ccc}%
0 & = & ...,0,0,0,0\\
&  & \\
1 & = & ...,0,0,0,1
\end{array}
\right.
\]
will be called \textbf{zero} and \textbf{one}, respectively.

\bigskip\ 

\begin{remark}
Please note that the dyadic integers and the rational integers lie in the set
of generic dyadic integers and the set of generic rational integers, respectively.
\end{remark}

\bigskip\ 

\begin{definition}
A \textbf{generic integer} $e$ such that
\[
e_{i}=1
\]
for almost all $i$ will be said to be \textbf{negative}. A generic integer $e
$ not equal $0$ such that
\[
e_{i}=0
\]
for almost all $i$ will be said to be \textbf{positive}. The generic integer
$e$ will be said to be \textbf{non-negative} if $e$ is either $0$ or positive.
\end{definition}

\bigskip\ 

\begin{remark}
Please note that there are generic rational integers which are neither
positive nor negative nor non-negative.
\end{remark}

\bigskip\ 

Let $a$ and $b$ be generic dyadic integers. The \textbf{component-wise
exclusive} \textbf{\textquotedblleft or\textquotedblright} of $a$ and $b$ ,
written $a\dotplus b$ , is defined as:
\[
a\dotplus b=...\;,\;a_{2}\dotplus b_{2},\;a_{1}\dotplus b_{1},\;a_{0}\dotplus
b_{0}%
\]
where $a_{i}\dotplus b_{i}$ denotes the exclusive \textquotedblleft
or\textquotedblright\ of the $i$-th components of $a$ and $b$ .

The \textbf{component-wise logical \textquotedblleft and\textquotedblright} of
$a$ and $b$ , written $a\diamond b$, is defined as:
\[
a\diamond b=...\;,\;a_{2}\diamond b_{2},\;a_{1}\diamond b_{1},\;a_{0}\diamond
b_{0}%
\]
where $a_{i}\diamond b_{i}$ denotes the logical \textquotedblleft
and\textquotedblright\ of the $i$-th components of $a$ and $b$ .

The \textbf{component-wise complement} of $a$ , written $a^{\ast}$ , is
defined as:
\[
a^{\ast}=...\;,\;a_{2}^{\ast},\;a_{1}^{\ast},\;a_{0}^{\ast}%
\]
where $a_{i}^{\ast}$ denotes the complement of the $i$-th component of $a$ .

Let $\alpha$ be an element of the free Boolean ring $\mathbb{B}\left\langle
\mathbf{x}\right\rangle $. Then the \textbf{scalar product} of $\alpha$ and
$a$ , written $\alpha\diamond a$, is defined as
\[
\alpha\diamond a=...\;,\;\alpha\diamond a_{2},\;\alpha\diamond a_{1}%
,\;\alpha\diamond a_{0},
\]
where $\alpha\diamond a_{i}$ denotes the logical \textquotedblleft
and\textquotedblright\ of $\alpha$ and the $i$-th component of $a$ .

The \textbf{unit left shift} of $a$ , written $\mathcal{S}a$, is defined as:
\[
\mathcal{S}a=...,a_{2},a_{1},a_{0},0.
\]
Let
\[
...\;,e^{(2)},\;e^{(1)},\;e^{(0)}%
\]
be an infinite sequence of elements of $\mathbb{G}_{(2)}\left\langle
\mathbf{x}\right\rangle $. If for every $j\geq0$, there exists a non-negative
integer $n(j)$ such that
\[
e_{j}^{(i)}=e_{j}^{\left(  n(j)\right)  }\qquad\qquad\text{for }i\geq n(j)
\]
then the sequence is said to be \textbf{convergent}. Otherwise, it is said to
be \textbf{divergent}. If the above sequence is convergent, its \textbf{limit}%
, written
\[
\lim_{i\rightarrow\infty}e^{(i)}%
\]
is said to \textbf{exist}, and is defined as the generic dyadic integer
\[
\lim_{i\rightarrow\infty}e^{(i)}=...\;,\;e_{2}^{\left(  n(2)\right)  }%
,\;e_{1}^{\left(  n(1)\right)  },\;e_{0}^{\left(  n(0)\right)  }%
\]

\bigskip\ 

\begin{definition}
An \textbf{instantiation} is a mapping
\[
\Phi:\mathbf{x}\longrightarrow\mathbb{B}\left\langle \ \right\rangle \text{ ,}%
\]
where $\mathbb{B}\left\langle \ \right\rangle $ denotes the free Boolean ring
on the empty set of symbols. Hence, $\mathbb{B}\left\langle \ \right\rangle $
may be identified with the field of two elements $\mathbb{F}_{2}$. Since
$\mathbb{B}\left\langle \mathbf{x}\right\rangle $ is free on $\mathbf{x}$ ,
every instantiation uniquely and naturally extends to a Boolean ring
epimorphism
\[
\Phi:\mathbb{B}\left\langle \mathbf{x}\right\rangle \longrightarrow
\mathbb{B}\left\langle \ \right\rangle
\]
which again is called an \textbf{instantiation}. Moreover, each instantiation
uniquely extends to epimorphisms:
\[
\left\{
\begin{array}
[c]{c}%
\Phi:\mathbb{G}\left\langle \mathbf{x}\right\rangle \longrightarrow
\mathbb{Z}\\
\\
\Phi:\mathbb{G}_{(2)}\left\langle \mathbf{x}\right\rangle \longrightarrow
\mathbb{Z}_{(2)}%
\end{array}
\right.
\]
which are also called \textbf{instantiations}.
\end{definition}

\bigskip\ 

The following will be helpful in proving theorems:

\bigskip\ 

\textbf{The Principle of Instantiation.}

\begin{description}
\item[a1)] Let $a$ and $b$ be elements of $\mathbb{B}\left\langle
\mathbf{x}\right\rangle $. If for every instantiation $\Phi$ ,
\[
\Phi(a)=\Phi(b)\;,
\]
then $a=b$.\bigskip

\item[a2)] Let $\Phi$ and $\Omega$ be instantiations. If for every element $a$
of $\mathbb{B}<\mathbf{x}>$,
\[
\Phi(a)=\Omega(a)\;,
\]
then $\Phi=\Omega$.\bigskip

\item[b1)] Let $a$ and $b$ be generic dyadic integers. If for every
instantiation $\Phi$,
\[
\Phi(a)=\Phi(b)\;,
\]
then $a=b$.\bigskip

\item[b2)] Let $\Phi$ and $\Omega$ be instantiations. If for every generic
dyadic integer $a$, ,
\[
\Phi(a)=\Omega(a)
\]
then $\Phi=\Omega$.
\end{description}

\bigskip

Finally, $\mathbb{G}\left\langle \mathbf{x}\right\rangle $ and $\mathbb{G}%
_{(2)}\left\langle \mathbf{x}\right\rangle $ can now be made into commutative
rings by defining two binary operations, \textbf{addition} \textquotedblleft%
$+$\textquotedblright\ and \textbf{multiplication} \textquotedblleft$\cdot
$\textquotedblright, as follows:

\bigskip\ 

\begin{definition}
[of addition \textquotedblleft+"]Let $a$ and $b$ be generic dyadic integers.
Let
\[
\left\{
\begin{array}
[c]{lll}%
c^{(0)} & = & a\\
&  & \\
carries^{(0)} & = & b
\end{array}
\right.
\]
and let
\[
\left\{
\begin{array}
[c]{lll}%
c^{(i+1)} & = & c^{(i)}\dotplus carries^{(i)}\\
&  & \\
carries^{(i+1)} & = & \mathcal{S}\left(  c^{(i)}\diamond carries^{(i)}\right)
\end{array}
\right.
\]
Then the sequences $c^{(i)}$ and $carries^{(i)}$ are convergent. The generic
dyadic integer $a+b$ is defined as:
\[
a+b=\lim_{i\rightarrow\infty}c^{(i)}%
\]
It can also be shown that
\[
\lim_{i\rightarrow\infty}carries^{(i)}=0
\]
Moreover, if $a$ and $b$ are generic rational integers, then $a+b$ is also a
generic rational integer.
\end{definition}

\bigskip

\begin{definition}
[of multiplication]Let $a$ and $b$ be generic dyadic integers. The
\textbf{product} of $a$ and $b$ , written $a\cdot b$, is defined as:
\[
a\cdot b=\sum_{i=0}^{\infty}b_{i}\diamond\left(  \mathcal{S}^{i}a\right)
=\sum_{i=0}^{\infty}a_{i}\diamond\left(  \mathcal{S}^{i}b\right)  \;,
\]
where \textquotedblleft$\sum_{i=0}^{\infty}$\textquotedblright\ denotes
\textquotedblleft$\lim_{j\rightarrow\infty}\sum_{i=0}^{j}$\textquotedblright,
and \textquotedblleft$\sum$\textquotedblright\ denotes a sum using the
operation \textquotedblleft$+$\textquotedblright\ defined above.
\end{definition}

\bigskip\ 

\begin{remark}
Generic integer multiplication \textquotedblleft$\cdot$\textquotedblright\ was
defined above in terms of the secondary operation of generic integer addition
\textquotedblleft$+$\textquotedblright. Please refer to the appendix for a
definition of generic integer multiplication \textquotedblleft$\cdot
$\textquotedblright\ in terms of more fundamental Boolean operations.
\end{remark}

\bigskip

Generic division \textquotedblleft$/$\textquotedblright\ will be defined in
the next section.

\bigskip

\section{Generic lopsided division and the Boolean factoring (BF) algorithm}

One of the objectives of this paper is to lift the algebraic structure (i.e.,
the fundamental binary operations) of the rational integers $\mathbb{Z}$ and
the dyadic integers $\mathbb{Z}_{(2)}$ to the generic level. \ In the previous
section, this was accomplished for all of the fundamental binary operations
but for the exception of division \textquotedblleft$/$\textquotedblright. \ In
this section, we complete this part of our research program by lifting the
lopsided division defined in section III to the generic level.

\bigskip

An immediate consequence of achieving his objective will be the creation of
the Boolean factoring (BF) algorithm, which transforms the problem of integer
factoring into the problem of solving a system of Boolean equations. \ This
system of Boolean equations is nothing more than the generic remainder arising
from generic lopsided division algorithm.

\bigskip\ 

\begin{definition}
Let $u$, $v$, $w$ be generic integers. The \textbf{first component-wise
symmetric function} of $u$, $v$, $w$ , written%
\[
Component\_Wise\_\sigma_{1}(u,v,w)
\]
is the generic integer whose $i$-th component is the first symmetric function
of the $i$-th components of $u$, $v$, $w$ , i.e., whose $i$-th component is
\[
\sigma_{1}(u,v,w)=u_{i}\dotplus v_{i}\dotplus w_{i}%
\]
The \textbf{second component-wise symmetric function} of $u$, $v$, $w$,
written
\[
Component\_Wise\_\sigma_{2}(u,v,w)
\]
is the generic integer whose $i$-th component is the second symmetric function
of the $i$-th components of $u$, $v$, $w$ , i.e., whose $i$-th component is
\[
\sigma_{2}(u_{i},v_{i},w_{i})=(u_{i}\diamond v_{i})\dotplus(v_{i}\diamond
w_{i})\dotplus(w_{i}\diamond u_{i})
\]

\end{definition}

\bigskip\ 

\begin{definition}
Let $u$ and $v$ be two generic rational integers. Let
\[
u\equiv v
\]
denote the following element of the free Boolean ring $\mathbb{B}%
<\mathbf{x}>$
\[
(u\equiv v)=\prod_{i=0}^{\infty}(1\dotplus u_{i}\dotplus v_{i}).
\]

\end{definition}

\bigskip\ 

\begin{remark}
Please note that, since $u$ and $v$ are generic integers, almost all terms in
the above product are $1$ .
\end{remark}

\bigskip\ 

In the theorem below, \textquotedblleft$\lfloor\;\rfloor$ \textquotedblright%
\ and \textquotedblleft$\lceil\;\rceil$\textquotedblright\ denote respectively
the \textbf{floor} and \textbf{ceiling} functions
\[%
\begin{tabular}
[c]{|c|}\hline
$\underset{\mathstrut}{\overset{\mathstrut}{%
\begin{array}
[c]{rcl}%
\lfloor\;\rfloor:\mathbb{R} & \longrightarrow & \mathbb{Z}\\
u & \mapsto & \max\left\{  k\in\mathbb{Z}:k\leq u\right\}
\end{array}
}}\text{ \ \ and \ \ }%
\begin{array}
[c]{rcl}%
\lceil\;\rceil:\mathbb{R} & \longrightarrow & \mathbb{Z}\\
u & \mapsto & \min\left\{  k\in\mathbb{Z}:k\geq u\right\}
\end{array}
$\\\hline
\end{tabular}
\ \text{ ,}%
\]
where $\mathbb{R}$ denotes the set of real numbers.

\bigskip\ 

\begin{theorem}
[Main](The Boolean Factoring Algorithm.) \ Let $N$ be a fixed positive
rational integer, and let
\[
N=\ldots\;,\;0,\;0,\;N_{\alpha-1},\;N_{\alpha-2},\;\ldots\;,\;N_{1},\;N_{0}%
\]
denote its radix 2 representation. Let
\[
\beta=\left\lfloor \left(  1+\alpha\right)  /2\right\rfloor \qquad,
\]
and let $x$ denote the positive odd generic rational integer
\[
x=\ldots\;,\;0,\;0,\;x_{\beta-1},\;x_{\beta-2},\;\ldots\;,\;x_{2},\;x_{1},\;1
\]
in $\mathbb{G}<\mathbf{x}>$ . (Hence, $x_{i}=0$ for $i\geq\beta$.) \ Let
\[
\left\{
\begin{array}
[c]{lll}%
c^{(0)} & = & N\\
&  & \\
borrows^{(0)} & = & 0
\end{array}
\right.
\]
and
\[
\left\{
\begin{array}
[c]{lll}%
c^{(i+1)} & = & Component\_Wise\_\sigma_{1}\left(  \;c^{(i)},\;borrows^{(i)}%
,\;c_{i}^{(i)}\diamond\mathcal{S}^{i}x\;\right) \\
&  & \\
borrows^{(i+1)} & = & Component\_Wise\_\sigma_{2}\left(  \;c^{(i)\ast
},\;borrows^{(i)},\;c_{i}^{(i)}\diamond\mathcal{S}^{i}x\;\right)
\end{array}
\right.
\]
where $c^{(i)\ast}$ denotes the component-wise complement of $c^{(i)}$ and
$\mathcal{S}$ denotes the unit left shift operator defined in section II of
this paper. \ Let
\[
\Gamma=\left\lceil \left(  1+\alpha\right)  /2\right\rceil \qquad.
\]
Finally, let $e_{k}$ denote the following element of $\mathbb{B}<\mathbf{x}>$
\[
e_{k}=\prod_{i=0}^{\infty}\left(  c_{i}^{(k)}\dotplus borrows_{i}%
^{(k)}\dotplus1\right)  =\left(  c^{(k)}\equiv borrows^{(k)}\right)
\]
for $k\geq0$ . Then $a$ has an odd rational integral factor of length
$\beta-j$ $(0\leq j\leq\beta-1)$ if and only if there exists an instantiation
$\Phi$ such that
\[
\left\{
\begin{array}
[c]{lllll}%
\Phi\left(  x_{\beta-p}\right)  & = & 0 &  & \text{for }1\leq p\leq j\\
&  &  &  & \\
\Phi\left(  x_{\beta-(j+1)}\right)  & = & 1 &  & \\
&  &  &  & \\
\Phi\left(  e_{\Gamma+1}\right)  & = & 1 &  &
\end{array}
\right.
\]
Moreover, if there exists such an instantiation $\Phi$ , then $N$ is the
product of the following two positive rational integers
\[
\Phi(x)
\]
and
\[
N/\Phi(x)=\Phi\left(  \;\ldots\;,\;0,\;c_{\Gamma+j-1}^{(\Gamma+j-1)}%
,\;c_{\Gamma+j-2}^{(\Gamma+j-2)},\;\ldots\;,\;c_{1}^{(1)},\;c_{0}%
^{(0)}\right)
\]

\end{theorem}

\bigskip

\begin{summary}
Thus, given an arbitrary positive integer $N$, the BF algorithm produces a
system of Boolean equations, namely
\[
c^{(\Gamma+1)}=borrows^{(\Gamma+1)}\text{ ,}%
\]
which we have expressed as the equality of two generic integers. \ Solving the
above system of Boolean equations, is equivalent to finding a satisfying set
of values for the Boolean variables $x_{1},x_{2},\ldots,x_{\beta-1}$ (i.e., an
instantiation $\Phi$) for the following single Boolean function%
\[
e_{\Gamma+1}=\prod_{i=0}^{\infty}\left(  c_{i}^{(\Gamma+1)}\dotplus
borrows_{i}^{(\Gamma+1)}\dotplus1\right)  \text{ ,}%
\]
i.e., finding a solution $\Phi$ such that%
\[
\Phi\left(  e_{\Gamma+1}\right)  =1\text{ .}%
\]
Each such satisfying set of values (i.e., each instantiation $\Phi$) produces
a rational integer divisor $\Phi(x)$ of the rational integer $N$, i.e.,%
\[
\Phi\left(  x\right)  /N\text{ .}%
\]

\end{summary}

\bigskip

\section{Examples of the application of the BF algorithm}

\bigskip

We\ now give a number of examples of integer factorization using the BF algorithm.

\bigskip

\begin{example}
Factoring $21$ with the Boolean factoring algorithm. The radix 2
representation of $21$ is:
\[
1\;0\;1\;0\;1
\]
Thus,
\[
\alpha=5,\;\beta=3,\;\Gamma=3\;.
\]%
\[%
\begin{tabular}
[c]{cccccccc}
&  &  & $c_{2}^{(2)}$ & $c_{1}^{(1)}$ & $c_{0}^{(0)}$ &  & \\
&  &  & $1\dotplus x_{2}$ & $x_{1}$ & $1$ & $a/x$ & \\\cline{2-6}\cline{6-6}%
$a=$ & $\ 1\ $ & $\ 0\ $ & $\ 1\ $ & $\ 0\ $ & $\ 1\ $ &
\multicolumn{1}{|c}{$x_{2}\;x_{1}\;1=x$} & \\
&  &  & $x_{2}$ & $x_{1}$ & $\ 1\ $ & $c_{0}^{(0)}\diamond\mathcal{S}^{0}x$ &
\\\cline{4-6}
&  &  & $1\dotplus x_{2}$ & $x_{1}$ & $\ 0\ $ & $c^{(1)}$ & \\
&  & $\ 0\ $ & $x_{1}$ &  &  & $borrows^{(1)}$ & \\
&  & $x_{1}x_{2}$ & $x_{1}$ & $x_{1}$ &  & $c_{1}^{(1)}\diamond\mathcal{S}%
^{1}x$ & \\\cline{3-5}
&  & $x_{1}x_{2}$ & $1\dotplus x_{2}$ & $\ 0\ $ &  & $c^{(2)}$ & \\
& $x_{1}x_{2}$ & $x_{1}$ &  &  &  & $borrows^{(2)}$ & \\
& $\ 0\ $ & $x_{1}\dotplus x_{1}x_{2}$ & $1\dotplus x_{2}$ &  &  &
$c_{2}^{(2)}\diamond\mathcal{S}^{2}x$ & \\\cline{2-4}
& $1\dotplus x_{1}x_{2}$ & $\ 0\ $ & $\ 0\ $ &  &  & $c^{(3)}$ & \\
& $x_{1}\dotplus x_{1}x_{2}$ &  &  &  &  & $borrows^{(3)}$ &
\end{tabular}
\ \ \ \ \ \
\]
Hence,
\[
e_{3}=1\dotplus x_{1}x_{2}\dotplus x_{1}\dotplus x_{1}x_{2}\dotplus1=x_{1}%
\]
and
\[
\left\{
\begin{array}
[c]{lll}%
\Phi(x_{2}) & = & 1\\
&  & \\
\Phi(x_{1}) & = & 1
\end{array}
\right.
\]
is a solution. Thus, $21$ is a product of the rational integers:
\[
\left\{
\begin{array}
[c]{lllll}%
\Phi\left(  x_{2},x_{1},1\right)  & = & 111\;\left(  base\;2\right)  & = &
7\;\left(  base\;10\right) \\
&  &  &  & \\
\Phi\left(  1\dotplus x_{2},x_{1},1\right)  & = & 011\;\left(  base\;2\right)
& = & 3\;\left(  base\;10\right)
\end{array}
\right.
\]

\end{example}

\bigskip\ 

\begin{example}
Factoring $77$ with the Boolean factoring algorithm.
\[
\alpha=7,\;\beta=4,\;\Gamma=4\qquad.
\]%
\[
\hspace{-0.75in}%
\begin{tabular}
[t]{cccccccc}
&  & \quad & $c_{3}^{(3)}$ & $c_{2}^{(2)}$ & $c_{1}^{(1)}$ & $c_{0}^{(0)}$ &
\\
&  & \quad & $1\dotplus x_{3}$ & $1\dotplus x_{2}$ & $x_{1}$ & $1$ &
$a/x$\\\cline{1-7}%
$1$ & $0$ & $0$ & $1$ & $1$ & $0$ & $1$ & \multicolumn{1}{|c}{$x_{3}%
,x_{2},x_{1},1=x$}\\
&  & \quad & $x_{3}$ & $x_{2}$ & $x_{1}$ & $1$ & $c_{0}^{(0)}\diamond
\mathcal{S}^{0}x$\\\cline{4-7}
&  & \quad & $1\dotplus x_{3}$ & $1\dotplus x_{2}$ & $x_{1}$ & $0$ & $c^{(1)}%
$\\
&  & $0$ &  & $x_{1}$ &  &  & $borrows^{(1)}$\\
&  & $x_{1}x_{3}$ & $x_{1}x_{2}$ & $x_{1}$ & $x_{1}$ &  & $c_{1}^{(1)}%
\diamond\mathcal{S}^{1}x$\\\cline{3-6}
&  & $x_{1}x_{3}$ & $1\dotplus x_{3}$ & $%
\begin{array}
[c]{c}%
1\dotplus x_{2}\\
\dotplus x_{1}x_{2}%
\end{array}
$ & $0$ &  & $c^{(2)}$\\
& $x_{1}x_{3}$ & $x_{1}x_{2}x_{3}$ & $x_{1}$ &  &  &  & $borrows^{(2)}$\\
& $x_{3}\dotplus x_{2}x_{3}$ & $0$ & $x_{1}\dotplus x_{1}x_{2}$ & $1\dotplus
x_{2}$ &  &  & $c_{2}^{(2)}\diamond\mathcal{S}^{2}x$\\\cline{2-5}
& $%
\begin{array}
[c]{c}%
x_{3}\dotplus x_{1}x_{3}\\
\dotplus x_{2}x_{3}%
\end{array}
$ & $x_{1}x_{3}\dotplus x_{1}x_{2}x_{3}$ & $1\dotplus x_{3}$ & $0$ &  &  &
$c^{(3)}$\\
$%
\begin{array}
[c]{c}%
x_{3}\dotplus x_{2}x_{3}\\
\dotplus x_{1}x_{2}x_{3}%
\end{array}
$ & $0$ & $x_{1}\dotplus x_{1}x_{2}x_{3}$ &  &  &  &  & $borrows^{(3)}$\\
$0$ & $x_{2}\dotplus x_{2}x_{3}$ & $x_{1}\dotplus x_{1}x_{3}$ & $1\dotplus
x_{3}$ &  &  &  & $c_{3}^{(3)}\diamond\mathcal{S}^{3}x$\\\cline{1-4}%
$%
\begin{array}
[c]{c}%
1\dotplus x_{3}\dotplus x_{2}x_{3}\\
\dotplus x_{1}x_{2}x_{3}%
\end{array}
$ & $x_{2}\dotplus x_{3}\dotplus x_{1}x_{3}$ & $0$ & $0$ &  &  &  & $c^{(4)}%
$\\
$x_{2}\dotplus x_{2}x_{3}$ & $x_{1}\dotplus x_{1}x_{3}$ &  &  &  &  &  &
$borrows^{(4)}$%
\end{tabular}
\ \ \ \ \ \ \ \ \ \ \ \ \
\]
Hence,
\begin{align*}
e_{4}  &  =\left[  \left(  1\dotplus x_{3}\dotplus x_{2}x_{3}\dotplus
x_{1}x_{2}x_{3}\right)  \dotplus\left(  x_{2}\dotplus x_{2}x_{3}\right)
\dotplus1\right]  \diamond\left[  \left(  x_{2}\dotplus x_{3}\dotplus
x_{2}x_{3}\right)  \dotplus\left(  x_{1}\dotplus x_{1}x_{3}\right)
\dotplus1\right] \\
& \\
&  =\left(  x_{2}\dotplus x_{3}\dotplus x_{1}x_{2}x_{3}\right)  \diamond
\left(  x_{1}\dotplus x_{2}\dotplus x_{3}\dotplus1\right) \\
& \\
&  =x_{1}\left(  x_{2}\dotplus x_{3}\right)
\end{align*}
Therefore, find an instantiation $\Phi$ such that
\[
\Phi\left(  x_{3}\right)  =1\text{ and }\Phi\left(  x_{1}(x_{2}\dotplus
x_{3})\right)  =1
\]
An algorithm for finding the solutions to Boolean equations of the above
\textbf{scarcely satisfiable} kind can be found in Section IV of this paper.
The solution $\Phi$ to these Boolean equations is:
\[
\Phi:\;%
\begin{array}
[c]{ccc}%
x_{1} & \longmapsto & 1\\
&  & \\
x_{2} & \longmapsto & 0\\
&  & \\
x_{3} & \longmapsto & 1
\end{array}
\]
Thus, $77$ is the product of the following positive rational integers:
\[
\Phi\left(  \;x_{3},\;x_{2},\;x_{1},\;1\;\right)  =1011\;\left(
base\;2\right)  =11\;\left(  base\;10\right)
\]
and
\[
\Phi\left(  \;1\dotplus x_{3},\;1\dotplus x_{2},\;x_{1},\;1\;\right)
=0111\;\left(  base\;2\right)  =7\;\left(  base\;10\right)
\]

\end{example}

\bigskip\ 

\begin{example}
Factoring $95$ with the Boolean factoring algorithm.
\[
\alpha=7,\;\beta=4,\;\Gamma=4
\]%
\[
\hspace{-0.75in}%
\begin{tabular}
[t]{cccccccc}
&  &  & $c_{3}^{(3)}$ & $c_{2}^{(2)}$ & $c_{1}^{(1)}$ & $c_{0}^{(0)}$ & \\
&  &  & $%
\begin{array}
[c]{c}%
1\dotplus x_{1}\\
\dotplus x_{2}\dotplus x_{3}%
\end{array}
$ & $1\dotplus x_{2}$ & $1\dotplus x_{1}$ & $\;1\;$ & $a/x$\\\cline{1-7}%
$1$ & $0$ & $1$ & $1$ & $1$ & $1$ & $1$ & \multicolumn{1}{|c}{$x_{3}%
,x_{2},x_{1},1$=$x$}\\
&  &  & $x_{3}$ & $x_{2}$ & $x_{1}$ & $1$ & $c_{0}^{(0)}\diamond
\mathcal{S}^{0}x$\\\cline{4-7}
&  &  & $1\dotplus x_{3}$ & $1\dotplus x_{2}$ & $1\dotplus x_{1}$ & $0$ &
$c^{(1)}$\\
&  & $0$ &  &  &  &  & $borrows^{(1)}$\\
&  & $x_{3}\dotplus x_{1}x_{3}$ & $x_{2}\dotplus x_{1}x_{2}$ & $0$ &
$1\dotplus x_{1}$ &  & $c_{1}^{(1)}\diamond\mathcal{S}^{1}x$\\\cline{3-6}
&  & $1\dotplus x_{3}\dotplus x_{1}x_{3}$ & $%
\begin{array}
[c]{c}%
1\dotplus x_{2}\\
\dotplus x_{3}\dotplus x_{1}x_{2}%
\end{array}
$ & $1\dotplus x_{2}$ & $0$ &  & $c^{(2)}$\\
& $0$ & $x_{2}x_{3}\dotplus x_{1}x_{2}x_{3}$ &  &  &  &  & $borrows^{(2)}$\\
& $x_{3}\dotplus x_{2}x_{3}$ & $0$ & $x_{1}\dotplus x_{1}x_{2}$ & $1\dotplus
x_{2}$ &  &  & $c_{2}^{(2)}\diamond\mathcal{S}^{2}x$\\\cline{2-5}
& $x_{3}\dotplus x_{2}x_{3}$ & $%
\begin{array}
[c]{c}%
1\dotplus x_{3}\dotplus x_{1}x_{3}\\
\dotplus x_{2}x_{3}\dotplus x_{1}x_{2}x_{3}%
\end{array}
$ & $%
\begin{array}
[c]{c}%
1\dotplus x_{1}\\
\dotplus x_{2}\dotplus x_{3}%
\end{array}
$ & $0$ &  &  & $c^{(3)}$\\
$x_{3}\dotplus x_{2}x_{3}$ & $x$$_{2}$$x$$_{3}\dotplus x_{1}x_{2}x_{3}$ &
$x_{1}x_{3}\dotplus x_{1}x_{2}x_{3}$ &  &  &  &  & $borrows^{(3)}$\\
$x$$_{1}$$x$$_{3}\dotplus x$$_{2}$$x$$_{3}$ & $x_{1}x_{2}\dotplus x_{2}x_{3}$
& $x_{1}x_{2}\dotplus x_{1}x_{3}$ & $%
\begin{array}
[c]{c}%
1\dotplus x_{1}\\
\dotplus x_{2}\dotplus x_{3}%
\end{array}
$ &  &  &  & $c_{3}^{(3)}\diamond\mathcal{S}^{3}x$\\\cline{1-4}%
$%
\begin{array}
[c]{c}%
1\dotplus x_{3}\\
\dotplus x_{1}x_{3}%
\end{array}
$ & $%
\begin{array}
[c]{c}%
x_{3}\dotplus x_{1}x_{2}\\
\dotplus x_{2}x_{3}\dotplus x_{1}x_{2}x_{3}%
\end{array}
$ & $%
\begin{array}
[c]{c}%
1\dotplus x_{3}\dotplus x_{1}x_{2}\\
\dotplus x_{1}x_{3}\dotplus x_{2}x_{3}%
\end{array}
$ & $0$ &  &  &  & $c^{(4)}$\\
$x_{1}x_{2}\dotplus x_{2}x_{3}$ & $x_{1}x_{3}\dotplus x_{1}x_{2}x_{3}$ &  &  &
&  &  & $borrows^{(4)}$%
\end{tabular}
\ \ \ \ \ \
\]

\end{example}

\bigskip\ 

The leftmost expression in $borrows^{(4)}$, i.e., $x_{1}x_{3}\dotplus
x_{1}x_{2}x_{3}$, is not shown in the above tableau because there is no room.

\bigskip\ 

Hence,
\begin{align*}
e_{4}  &  =\left(  0\dotplus x_{1}x_{3}\dotplus x_{1}x_{2}x_{3}\dotplus
1\right)  \diamond\left(  1\dotplus x_{3}\dotplus x_{1}x_{3}\dotplus
x_{1}x_{2}\dotplus x_{2}x_{3}\dotplus1\right) \\
&  \diamond\left(  x_{3}\dotplus x_{1}x_{2}\dotplus x_{2}x_{3}\dotplus
x_{1}x_{2}x_{3}\dotplus x_{1}x_{3}\dotplus x_{1}x_{2}x_{3}\dotplus1\right) \\
&  \diamond\left(  1\dotplus x_{3}\dotplus x_{1}x_{2}\dotplus x_{1}%
x_{3}\dotplus x_{2}x_{3}\dotplus0\dotplus1\right)  \diamond1\diamond1
\end{align*}
Therefore,
\begin{align*}
e_{4}  &  =\left(  x_{1}x_{3}\dotplus x_{1}x_{2}x_{3}\dotplus1\right)
\diamond\left(  x_{3}\dotplus x_{1}x_{3}\dotplus x_{1}x_{2}\dotplus x_{2}%
x_{3}\right) \\
&  \diamond\left(  x_{3}\dotplus x_{1}x_{2}\dotplus x_{1}x_{3}\dotplus
x_{2}x_{3}\dotplus1\right)  \diamond\left(  x_{3}\dotplus x_{1}x_{2}\dotplus
x_{1}x_{3}\dotplus x_{2}x_{3}\right)
\end{align*}
Hence,
\[
e_{4}=0
\]
Thus, there is no factor of $95$ of length $4$ .

\bigskip\ 

So we set $x_{3}=0$ , and continue.

\bigskip\ %

\[
\hspace{-0.5in}%
\begin{tabular}
[t]{ccccccccc}
&  &  &  & $c_{3}^{(3)}$ & $c_{2}^{(2)}$ & $c_{1}^{(1)}$ & $c_{0}^{(0)}$ & \\
&  &  &  & $1\dotplus x_{1}\dotplus x_{2}$ & $1\dotplus x_{2}$ & $1\dotplus
x_{1}$ & $\;1\;$ & $a/x$\\\cline{1-8}
& $1$ & $0$ & $1$ & $1$ & $1$ & $1$ & $1$ & \multicolumn{1}{|c}{$0,x_{2}%
,x_{1},1=x$}\\
&  &  &  & $0$ & $x_{2}$ & $x_{1}$ & $1$ & $c_{0}^{(0)}\diamond\mathcal{S}%
^{0}x$\\\cline{5-8}
&  &  &  & $1$ & $1\dotplus x_{2}$ & $1\dotplus x_{1}$ & $0$ & $c^{(1)}$\\
&  &  & $0$ &  &  &  &  & $borrows^{(1)}$\\
&  &  & $0$ & $x_{2}\dotplus x_{1}x_{2}$ & $0$ & $1\dotplus x_{1}$ &  &
$c_{1}^{(1)}\diamond\mathcal{S}^{1}x$\\\cline{4-7}
&  &  & $1$ & $1\dotplus x_{2}\dotplus x_{1}x_{2}$ & $1\dotplus x_{2}$ & $0$ &
& $c^{(2)}$\\
&  & $0$ & $0$ &  &  &  &  & $borrows^{(2)}$\\
&  & $0$ & $0$ & $x_{1}\dotplus x_{1}x_{2}$ & $1\dotplus x_{2}$ &  &  &
$c_{2}^{(2)}\diamond\mathcal{S}^{2}x$\\\cline{3-6}
&  & $0$ & $1$ & $1\dotplus x_{1}\dotplus x_{2}$ & $0$ &  &  & $c^{(3)}$\\
& $0$ & $0$ & $0$ &  &  &  &  & $borrows^{(3)}$\\
& $0$ & $x_{1}x_{2}$ & $x_{1}x_{2}$ & $1\dotplus x_{1}\dotplus x_{2}$ &  &  &
& $c_{3}^{(3)}\diamond\mathcal{S}^{3}x$\\\cline{1-5}
& $1$ & $x_{1}x_{2}$ & $1\dotplus x_{1}x_{2}$ & $0$ &  &  &  & $c^{(4)}$\\
$0$ & $x_{1}x_{2}$ & $0$ &  &  &  &  &  & $borrows^{(4)}$\\
$0$ & $x_{2}\dotplus x_{1}x_{2}$ & $x$$_{1}\dotplus x$$_{1}$$x$$_{2}$ &
$1\dotplus x_{1}x_{2}$ &  &  &  &  & $c_{4}^{(4)}\diamond\mathcal{S}^{4}%
x$\\\cline{1-4}
& $1\dotplus x_{2}$ & $x_{1}$ & $0$ &  &  &  &  & $c^{(5)}$\\
$0$ & $x_{1}\dotplus x_{1}x_{2}$ &  &  &  &  &  &  & $borrows^{(5)}$%
\end{tabular}
\ \ \ \ \ \
\]

\bigskip\ 

Therefore,
\begin{align*}
e_{5}  &  =\left(  1\dotplus x_{2}\dotplus x_{1}\dotplus x_{1}x_{2}%
\dotplus1\right)  \diamond\left(  x_{1}\dotplus0\dotplus1\right)
\diamond1\diamond1\diamond\;\cdots\;\diamond1\\
& \\
&  =\left(  x_{1}\dotplus x_{2}\dotplus x_{1}x_{2}\right)  \diamond\left(
1\dotplus x_{1}\right)  =x_{2}\left(  1\dotplus x_{1}\right)
\end{align*}
Thus, a solution $\Phi$ is:
\[
\Phi:\;%
\begin{array}
[c]{c}%
x_{1}\longmapsto0\\
\\
x_{2}\longmapsto1\\
\\
x_{3}\longmapsto0
\end{array}
\]
Hence, $95$ is the product of:
\[
\Phi\left(  \;x_{3},\;x_{2},\;x_{1},\;1\;\right)  =0101\;\left(
base\;2\right)  =5\;\left(  base\;10\right)
\]
and
\[
\Phi\left(  \;1\dotplus x_{1}x_{2},\;1\dotplus x_{1}\dotplus x_{2},\;1\dotplus
x_{2},\;1\dotplus x_{1},\;1\;\right)  =10011\;\left(  base\;2\right)
=19\;\left(  base\;10\right)
\]

\bigskip

\section{The multiplicative Boolean factoring (MBF) algorithm}

\bigskip

The BF algorithm defined in the previous section is based on generic lopsided
division. \ But there is no need to perform generic lopsided division each
time one factors an integer. \ One need only pre-compute the generic inverse
$x^{-1}$ of a judiciously chosen odd generic integer $x$, and then use the
pre-computed inverse $x^{-1}$ over and over again to factor arbitrarily chosen integers.

\bigskip

Let $x$ denote the following particular generic integer%
\[
x=\ldots,x_{3},x_{2},x_{1},1
\]
and let $x^{-1}$ be the corresponding generic inverse, which can be computed
with generic lopsided division.

\bigskip

Let $x$ be pre-computed. \ Then for each chosen positive integer $N$ to be
factored, one can find the appropriate system of Boolean equations to be
solved to factor $N$ simply by computing the generic product%
\[
N\cdot x^{-1}\text{ .}%
\]
The resulting algorithm is called the \textbf{Multiplicative Boolean factoring
(MBF) algorithm}.

\bigskip

We leave the remaining details to the reader.

\bigskip

\section{A method for solving scarcely satisfiable Boolean equations}

\bigskip

In this section, we outline a general framework for solving the system of
Boolean equations produced by the BF and MBF algorithms. \ This framework was
later used by Gamal Abdali to create a LISP implementation of the BF
algorithm. \ Sumeet Bagde then extended these methods by using binary decision
diagrams (BDDs)\cite{Bryant2} to create a Mathematica program that also
implemented the BF algorithm.

\bigskip

Both the LISP and Mathematica programs were used to factor many integers.
\ The runtime statistics clearly indicated that the BF algorithm runs in
exponential time, and hence, is not competitive with the best classical
factoring algorithms. \ For an algebraic proof as to why this is the case, we
refer the reader to the topdown overview given in the next section of this paper.

\bigskip

The main theorem, found in section IV, reduces the task of factoring a fixed
positive rational integer $N$ to the task of finding a solution to a Boolean
equation of the form
\[
e=1
\]
where $e\in\mathbb{B}<\mathbf{x}>$. Each solution of this equation corresponds
to a divisor of $N$ .

On the other hand, we are interested in factoring large integers which only
have a \textbf{small} number of divisors. It follows that the corresponding
equation
\[
e=1
\]
has only a \textbf{small} number of solutions. We now use this idea to develop
a method for solving Boolean equations that each have only a small number of solutions.

\bigskip\ 

\begin{definition}
Let $e\in\mathbb{B}\left\langle \mathbf{x}\right\rangle $. A \textbf{solution}
of the Boolean equation
\[
e=1
\]
is an instantiation $\Phi$ such that
\[
\Phi(e)=1\qquad.
\]

\end{definition}

\bigskip\ 

\begin{definition}
Let \textbf{less than}, written \textquotedblleft$<$\textquotedblright, denote
the linear ordering on the free basis elements
\[
\mathbf{x=}\left\{  \;\ldots\;,\;x_{2},\;x_{1},\;x_{0}\;\right\}
\]
defined by
\[
x_{i}<x_{j}\;\;\text{if}\;\;i<j
\]
A \textbf{term} in $\mathbb{B}\left\langle \mathbf{x}\right\rangle $ is a
finite product of distinct free basis elements which appear in the product
from left-to-right in ascending order according to the relation
\textquotedblleft$<$\textquotedblright\ . The element $1$ of $\mathbb{B}%
\left\langle \mathbf{x}\right\rangle $ is represented as the term which is the
empty product of free basis elements. Let \textquotedblleft$<$%
\textquotedblright\ also denote the lexicographic linear ordering induced on
the set of terms by the linear ordering \textquotedblleft$<$\textquotedblright%
\ on $\mathbf{x}$. (Please note that $1$ is the smallest term.) A
\textbf{canonical expression} is a sum of distinct terms which appear in the
sum from left-to-right in ascending order according to the linear ordering
\textquotedblleft$<$\textquotedblright\ . (Please note that $0$ is represented
by the empty sum of terms.)
\end{definition}

\bigskip\ 

\noindent\textbf{Observation.} Every element of $\mathbb{B}\left\langle
\mathbf{x}\right\rangle $ is uniquely representable as a canonical expression.

\bigskip\ 

Finally, we are in a position to define what is meant by a Boolean equation
having only a small number of solutions.

\bigskip\ 

\begin{definition}
Let $e\in\mathbb{B}\left\langle \mathbf{x}\right\rangle $ . The Boolean
equation
\[
e=1
\]
is said to be \textbf{scarcely satisfiable} if the number of its solutions
$\Phi$ is a non-zero number which is less than the number of distinct free
basis elements $x_{i}$ appearing in the canonical expression for $e$ .
\end{definition}

\bigskip\ 

Next, we observe that the solutions of
\[
e=1
\]
are in 1-1 correspondence with minterms in the minterm expansion of $e$ . It
follows that a scarcely satisfiable Boolean equation
\[
e=1
\]
is one in which there are only a small number of minterms appearing in the
minterm expansion of $e$ . Thus, we have the following proposition.

\bigskip\ 

\begin{proposition}
Let
\[
e=1
\]
be a scarcely satisfiable Boolean equation. Then for all but a small number of
free basis elements $x_{i}$ appearing in the canonical expression for $e$
either
\[
x_{i}e=e
\]
or
\[
x_{i}^{\ast}e=e
\]

\end{proposition}

\bigskip\ 

\begin{remark}
Please note that
\end{remark}

\begin{description}
\item[1)] $x_{i}e=0\Longleftrightarrow x_{i}^{\ast}e=e\Longrightarrow
\Phi(x_{i})=0$ for all solutions $\Phi$ of $e=1$ .

\item[2)] $x_{i}^{\ast}e=0\Longleftrightarrow x_{i}e=e\Longrightarrow
\Phi(x_{i})=1$ for all solutions $\Phi$ of $e=1$ .
\end{description}

\bigskip\ 

This leads to the following:

\bigskip\ 

\noindent\textbf{Algorithm for finding a solution }$\Phi$ \textbf{to a
scarcely satisfiable Boolean equation}
\[
e=1
\]

\begin{description}
\item[Step 1] For each free basis element $x_{i}$ not appearing in the
canonical expression for $e$, set $\Phi(x_{i})$ arbitrarily equal to $0$ or $1
$ .

\item[Step 2] For each free basis element $x_{i}$ such that
\[
x_{i}^{\ast}e=0
\]
set
\[
\Phi(x_{i})=1
\]

\item[Step 3] For each free basis element $x_{i}$ such that
\[
x_{i}e=0
\]
set
\[
\Phi(x_{i})=0
\]

\item[Step 4] Let $\mu$ denote the number of free basis elements not
determined in Steps 1 through 4. Exhaustively try each of the possible
$2^{\mu}$ assignments of $\Phi$ for these basis elements. Since $e$ is
scarcely satisfiable, the number of possibilities $2^{\mu}$ is small.
\end{description}

\bigskip

\section{A topdown overview: The "big picture"}

\bigskip

After the careful microscopic analysis of arithmetic complexity given in the
previous sections of this paper, the "big picture" emerges. \ We now step
back, and take a discerning macroscopic look at what has been found. \ 

\bigskip

We begin by defining two rings $\mathbb{B}^{\oplus}$ and $\mathbb{Z}^{\oplus}%
$, the former of characteristic $2$, the latter of characteristic $0$.

\bigskip

\noindent\textsc{Construction of the first ring:}

\bigskip

Let $\mathbb{B}\left\langle \underline{\underline{x}}\right\rangle $ be the
free Boolean ring on the set
\[
\underline{\underline{x}}=\left\{  x_{0},x_{1},\ldots,x_{n-1}\right\}
\]
of $n$ symbols. \ Since the ring $\mathbb{B}\left\langle
\underline{\underline{x}}\right\rangle $ is both semisimple and a principle
ideal domain, it decomposes into the direct sum%
\[
\mathbb{B}\left\langle \underline{\underline{x}}\right\rangle =%
{\displaystyle\bigoplus\limits_{\alpha=0}^{2^{n}-1}}
\left(  m_{\alpha}\right)  \simeq%
{\displaystyle\bigoplus\limits_{\alpha=0}^{2^{n}-1}}
\mathbb{F}_{2}%
\]
of principal, minimal ideals, where the ideal generators $m_{\alpha}$, called
\textbf{minterms}, form a complete set of orthogonal idempotents.
\ Consequently, each element $e\in\mathbb{B}\left\langle
\underline{\underline{x}}\right\rangle $ can be uniquely written in the form%
\[
e=%
{\displaystyle\sum\limits_{\alpha=0}^{2^{n}-1}}
c_{\alpha}m_{a}\text{ ,}%
\]
where\ each coefficient $c_{\alpha}$ is an element of the finite field of two
elements $\mathbb{F}_{2}=\left\{  0,1\right\}  $.

\bigskip

We now define the \textbf{ring }$\mathbb{B}^{\oplus}$\textbf{ of concurrent
Boolean functions} to be the Boolean ring formed by the ring direct sum%
\[
\mathbb{B}^{\oplus}=%
{\displaystyle\bigoplus\limits_{j=0}^{\infty}}
\mathbb{B}\left\langle \underline{\underline{x}}\right\rangle
\]
with diagonal multiplication. \ 

\bigskip

\noindent\textsc{Construction of the second ring:}

\bigskip

Let $\mathbb{Z}$ be the ring of rational integers. We define the \textbf{ring
}$\mathbb{Z}^{\oplus}$ \textbf{of concurrent integers} as the direct sum
\[
\mathbb{Z}^{\oplus}=\mathbb{Z}^{\oplus2^{n}}=%
{\displaystyle\bigoplus\limits_{\alpha=0}^{2^{n}-1}}
\mathbb{Z}%
\]
with multiplication defined diagonally. \ Thus each element of $\mathbb{Z}%
^{\oplus}$ can be uniquely written in the form
\[%
{\textstyle\bigoplus\limits_{\alpha=0}^{2^{n}-1}}
b_{\alpha}=%
{\textstyle\bigoplus\limits_{\alpha=0}^{2^{n}-1}}
{\displaystyle\sum\limits_{j=0}^{\infty}}
b_{j\alpha}2^{j}\text{ ,}%
\]
where $b_{\alpha}\in\mathbb{Z}$, and where
\[
b_{\alpha}=%
{\displaystyle\sum\limits_{j=0}^{\infty}}
b_{j\alpha}2^{j}%
\]
is the binary expansion of the integer $b_{\alpha}$.

\bigskip

\noindent\textsc{The identification:}

\bigskip

We now have two rings, $\mathbb{B}^{\oplus}$ of characteristic $2$, and
$\mathbb{Z}^{\oplus}$ of characteristic $0$. \ Our next step is to identify
these two rings \underline{as} \underline{sets} via the bijection defined by%

\[%
\begin{array}
[c]{ccc}%
\mathbb{Z}^{\oplus} &
\begin{array}
[c]{c}%
\overset{\Upsilon}{\longrightarrow}\\
\underset{\Upsilon^{-1}}{\longleftarrow}%
\end{array}
& \mathbb{B}^{\oplus}\\
&  & \\%
{\textstyle\bigoplus\limits_{\alpha=0}^{2^{n}-1}}
{\displaystyle\sum\limits_{j=0}^{\infty}}
b_{j\alpha}2^{j} &
\begin{array}
[c]{c}%
\overset{\Upsilon}{\longmapsto}\\
\underset{\Upsilon^{-1}}{\longleftarrow}%
\end{array}
&
{\textstyle\bigoplus\limits_{j=0}^{\infty}}
{\displaystyle\sum\limits_{\alpha=0}^{2^{n}-1}}
b_{j\alpha}m_{\alpha}\text{ ,}%
\end{array}
\]
where $b_{j\alpha}$ on the left is an integer $0$ or $1$ in $\mathbb{Z}$, and
on the right an element $0$ or $1$ of the finite field of two elements
$\mathbb{F}_{2}$. \ The result of this identification is called the
\textbf{ring of generic integers}, and denoted by $\mathbb{G}\left\langle
\underline{\underline{x}}\right\rangle $. \ This object $\mathbb{G}%
\left\langle \underline{\underline{x}}\right\rangle $ is a set with two
distinct ring structure, i.e., a \textbf{bi-ring} $\mathbb{G}\left\langle
\underline{\underline{x}}\right\rangle ,+,\cdot,\dotplus,\diamond$. \ The
bi-ring $\mathbb{G}\left\langle \underline{\underline{x}}\right\rangle $ is
simultaneously of characteristic $0$ and of characteristic $2$.

\bigskip

\begin{remark}
Please take care to note that both $\Upsilon$ and its inverse $\Upsilon^{-1}$
are bijections of sets, and not ring isomorphisms.
\end{remark}

\bigskip

\bigskip

In like manner, the bi-ring of generic dyadic integers $\mathbb{G}%
_{(2)}\left\langle \underline{\underline{x}}\right\rangle $ can be defined as
follows: \ 

\bigskip

Let $\mathbb{Z}_{(2)}$ denote the ring of dyadic integers. \ We define the
(characteristic $0$) \textbf{ring }$\mathbb{Z}_{(2)}^{\times}$ of
\textbf{concurrent dyadic integers} as the direct product%
\[
\mathbf{\ }\mathbb{Z}_{(2)}^{\times}=\overset{2^{n}-1}{\underset{\alpha
=0}{\times}}\mathbb{Z}_{(2)}\text{ }%
\]
with diagonal multiplication. \ 

\bigskip

In turn, the (characteristic $2$) \textbf{dyadic ring} $\mathbb{B}%
_{(2)}^{\times}$ \textbf{of concurrent Boolean functions} is defined as the
direct product%
\[
\mathbb{B}_{(2)}^{\times}=\overset{\infty}{\underset{j=0}{\times}}%
\mathbb{B}\left\langle \underline{\underline{x}}\right\rangle
\]
with diagonal multiplication. \ 

\bigskip

The definition of the bijection $\Upsilon_{(2)}$, i.e.,
\[%
\begin{array}
[c]{ccc}%
\mathbb{Z}_{(2)}^{\times} &
\begin{array}
[c]{c}%
\overset{\Upsilon_{(2)}}{\longrightarrow}\\
\underset{\Upsilon_{(2)}^{-1}}{\longleftarrow}%
\end{array}
& \mathbb{B}_{(2)}^{\times}%
\end{array}
\text{ .}%
\]
is similar to that of $\Upsilon$, and is left to the reader.

\bigskip

Finally, the \textbf{bi-ring of generic dyadic integers} $\mathbb{G}%
_{(2)}\left\langle \underline{\underline{x}}\right\rangle ,+,\cdot
,\dotplus,\diamond$ is defined by using the above bijection $\Upsilon_{(2)}$
to identify $\mathbb{B}_{(2)}^{\times}$ and $\mathbb{Z}_{(2)}^{\times}$ as sets.

\bigskip

One immediate consequence of the above "big picture" is that our use of the
phrase "algebraic parallel processing" is not unwarranted. \ For each
arithmetic operation in $\mathbb{G}\left\langle \underline{\underline{x}%
}\right\rangle ,+,\cdot$ \ (in $\mathbb{G}_{(2)}\left\langle
\underline{\underline{x}}\right\rangle ,+,\cdot$) is the same as $2^{n}-1$
arithmetic operations in $\mathbb{G}\left\langle \underline{\underline{x}%
}\right\rangle ,\dotplus,\diamond$ \ (in $\mathbb{G}_{(2)}\left\langle
\underline{\underline{x}}\right\rangle ,\dotplus,\diamond$).

\bigskip

It also follows that we have achieved one of the research objectives mentioned
in the introduction, namely, the development of a procedure for decomposing
arithmetic operations into fundamental Boolean operations. \ In other words,
we have developed a procedure for decomposing the arithmetic operations of the
ring $\mathbb{G}\left\langle \underline{\underline{x}}\right\rangle ,+,\cdot$
\ (of the ring $\mathbb{G}_{(2)}\left\langle \underline{\underline{x}%
}\right\rangle ,+,\cdot$) into the elementary operations of the ring
$\mathbb{G}\left\langle \underline{\underline{x}}\right\rangle ,\dotplus
,\diamond$ \ (of the ring $\mathbb{G}_{(2)}\left\langle
\underline{\underline{x}}\right\rangle ,\dotplus,\diamond$).

\vspace{0.5in}

In closing this section, we give below a summary of the Boolean decompositions
of the elementary arithmetic operations for addition "$+$", negation "$-$",
subtraction "$-$", and lopsided division "$/$". \ The Boolean decomposition of
multiplication "$\cdot$" can be found in the appendix. \ Please note that all
of the operations $+$, $-$, $/$ can be viewed as fixed points of either the
function $A$ or the function $D$. \ The definitions of the functions $A$ and
$D$ can be found below.

\begin{proposition}
Let $A$ and $P$ be the functions defined by
\[
\hspace{-0.75in}%
\begin{array}
[c]{ccc}%
A:\mathbb{G}_{(2)}\left\langle \underline{\underline{x}}\right\rangle
\times\mathbb{G}_{(2)}\left\langle \underline{\underline{x}}\right\rangle  &
\longrightarrow & \mathbb{G}_{(2)}\left\langle \underline{\underline{x}%
}\right\rangle \times\mathbb{G}_{(2)}\left\langle \underline{\underline{x}%
}\right\rangle \\
\left(  u,v\right)  & \longmapsto & \left(  a\dotplus b,\mathcal{S}\left(
a\diamond b\right)  \right)
\end{array}
\text{ \ \ and \ \ }%
\begin{array}
[c]{ccc}%
P:\mathbb{G}_{(2)}\left\langle \underline{\underline{x}}\right\rangle
\times\mathbb{G}_{(2)}\left\langle \underline{\underline{x}}\right\rangle  &
\longrightarrow & \mathbb{G}_{(2)}\left\langle \underline{\underline{x}%
}\right\rangle \\
\left(  u,v\right)  & \longmapsto & u
\end{array}
\text{ .}%
\]
Then

\begin{itemize}
\item[\textsc{Addition} "$+$"]
\[
a+b=P\lim_{k\rightarrow\infty}A^{k}\left(  a,b\right)  \text{ ,}%
\]
and hence, $\left(  a+b,0\right)  $ is a fixed point of $A.$

\item[\textsc{Negation} "$-$"]
\[
-a=P\lim_{k\rightarrow\infty}A^{k}\left(  a^{\ast},1\right)  \text{ ,}%
\]
and hence, $\left(  -a,0\right)  $ is a fixed point of the function $A$.

\item[\textsc{Subtraction} "$-$"]
\[
a-b=P\lim_{k\rightarrow\infty}A^{k}\left(  a\dotplus b^{\ast}\dotplus
1,\mathcal{S}\left(  a\diamond b^{\ast}\dotplus a_{0}\dotplus b_{0}^{\ast
}\right)  \right)  \text{ ,}%
\]
and hence, $\left(  a-b,0\right)  $ is a fixed point of the function $A$.
\end{itemize}

\bigskip

Let $\sigma_{1}$ and $\sigma_{2}$ are the component-wise symmetric functions
defined in section IV. \ Let $D$ and $P^{\prime}$ be the functions defined by%
\[%
\begin{array}
[c]{ccc}%
D:\mathbb{G}_{(2)}\left\langle \underline{\underline{x}}\right\rangle
^{3}\times\mathbb{Z} & \longrightarrow & \mathbb{G}_{(2)}\left\langle
\underline{\underline{x}}\right\rangle ^{3}\times\mathbb{Z}\\
\left(  u,v,w,\ell\right)  & \longmapsto & \left(  \sigma_{1}\left(
u,v,u_{\ell}\diamond w\right)  ,\sigma_{2}\left(  u^{\ast},v,u_{\ell}\diamond
w\right)  ,\mathcal{S}\left(  w\right)  ,\ell+1\right)
\end{array}
,
\]
and%
\[%
\begin{array}
[c]{ccc}%
P^{\prime}:\mathbb{G}_{(2)}\left\langle \underline{\underline{x}}\right\rangle
^{3}\times\mathbb{Z} & \longrightarrow & \mathbb{G}_{(2)}\left\langle
\underline{\underline{x}}\right\rangle \\
\left(  u,v,w,\ell\right)  & \longmapsto & u
\end{array}
\text{ .}%
\]

\end{proposition}

Then\bigskip

\textsc{Lopsided Division} "$/$"%

\[
a/b=P^{\prime}\lim_{k\rightarrow\infty}D^{k}\left(  a,0,b,0\right)  \text{ ,}%
\]
\textit{and hence, }$\left(  a/b,0,0,\infty\right)  $\textit{ is a fixed point
of }$D$\textit{, where }$a$\textit{ and }$b$\textit{ are odd positive generic
dyadic integers. }

\bigskip

\noindent\textbf{Question:} \ \textit{Is it possible to remove the
characteristic }$0$\textit{ counter }$\ell\longmapsto\ell+1$\textit{ from the
above lopsided division algorithm?}

\bigskip

\section{Conclusions and open questions}

\bigskip

The BF and MBF algorithms described in this paper are far from competitive
with current classical factoring algorithms. \ It is hoped that the two
Boolean factoring algorithms and the generic framework described within this
paper will become natural stepping stones for creating faster and more
competitive algorithms, and perhaps lead to a new highly competitive quantum
integer factorization algorithm.

\bigskip

Before closing this section, a word should be said about the natural question
of what is the relationship between the satisfiability problem SAT and the BF
and MBF algorithms. \ SAT is NP-complete. \ Is the task of solving the class
of systems of Boolean equations produced by the BF and MBF algorithms (when
reduced to a decision problem) NP-complete? \ Is it \#P-hard? \ The answers to
these questions are not known at this time.

\bigskip

Certainly the BF and MBF algorithms culminate in satisfiability problems.
\ But not all satisfiability problems are NP-complete. \ Not all are \#P hard.
\ For example, 2-SAT is not NP-complete, but 3-SAT is.

\bigskip

It should be noted that the execution of every algorithm on a digital computer
is ultimately reduced by\ a compiler/assembler to the execution of a Boolean
algorithm. \ This was one of the primary motivating factors for writing this
paper. \ Certainly, it is clear that not every Boolean algorithm executed on a
digital computer corresponds to an NP-complete or a \#P hard problem.

\bigskip

\section{Appendix. Generic integer multiplication defined in terms of
fundamental Boolean operations}

In section II, generic integer multiplication \textquotedblleft$\cdot
$\textquotedblright\ was defined in terms of the secondary operation of
generic integer addition \textquotedblleft$+$\textquotedblright. Sketched
below is a definition of generic integer multiplication \textquotedblleft%
$\cdot$\textquotedblright\ in terms of more fundamental Boolean operations.

\bigskip\ 

\begin{definition}
Let
\[
a=\ldots\;,\;a_{2},\;a_{1},\;a_{0}%
\]
be a non-negative (i.e., positive or zero) generic rational integer. Then the
$i$-\textbf{th} \textbf{symmetric function}
\[
\sigma_{i}(a)
\]
of $a$ is defined as
\[
\sigma_{i}(a)=\sum_{j(1)<j(2)<\;\ldots\;<j(i)}a_{j(1)}\diamond a_{j(2)}%
\diamond\;\ldots\;\diamond a_{j(i)}%
\]
where
\[
\sum_{j(1)<j(2)<\;\ldots\;<j(i)}%
\]
denotes a sum with respect to the operation \textquotedblleft$\dotplus
$\textquotedblright\ in $B<\mathbf{x}>$ over the indices $j(1),\;j(2),\;\ldots
\;,\;j(i)$ subject to the condition
\[
j(1)<j(2)<\;\ldots\;<j(i)\;.
\]
This function is well-defined since
\[
a_{k}=0
\]
for all but finitely many $k$.
\end{definition}

\bigskip\ 

\begin{definition}
Let
\[
\ldots\;,\;a^{(2)},\;a^{(1)},\;a^{(0)}%
\]
denote an infinite sequence of generic dyadic integers such that
\[
\lim_{i\rightarrow\infty}a^{(i)}=0
\]
Then, since the limit of the sequence is zero, the $j$-th components
\[
\ldots\;,\;a_{j}^{(2)},\;a_{j}^{(1)},\;a_{j}^{(0)}%
\]
of the above elements of the sequence form a generic rational integer. The
\textbf{component-wise }$i$\textbf{-symmetric function}
\[
Component\_Wise\_\sigma_{i}\left(  \ldots\;,\;a^{(2)},\;a^{(1)},\;a^{(0)}%
\right)
\]
of the sequence
\[
\ldots\;,\;a^{(2)},\;a^{(1)},\;a^{(0)}%
\]
is defined as the generic dyadic integer whose $j$-th component is
\[
\sigma_{i}\left(  \ldots\;,\;a_{j}^{(2)},\;a_{j}^{(1)},\;a_{j}^{(0)}\right)
\]

\end{definition}

\bigskip\ 

\noindent\textbf{Observation.} Let
\[
\ldots\;,\;a^{(2)},\;a^{(1)},\;a^{(0)}%
\]
be an arbitrary sequence of generic dyadic integers, and let \textquotedblleft%
$\sum$\textquotedblright\ denote a summation with respect to the operation
\textquotedblleft$+$\textquotedblright. Then

\begin{description}
\item[1)] $\sum_{i=0}^{\infty}a^{(i)}$ is convergent, hence exists, if and
only if $\lim_{i\rightarrow\infty}a^{(i)}=0$.

\item[2)] $\lim_{i\rightarrow\infty}\mathcal{S}^{i}a^{(i)}=0$, and hence,
$\sum_{i=0}^{\infty}\mathcal{S}^{i}a^{(i)}$ is convergent and well-defined.
\end{description}

\bigskip\ 

\begin{proposition}
Let
\[
\ldots\;,\;c^{(2)},\;c^{(1)},\;c^{(0)}%
\]
be a sequence of generic dyadic integers. Define $c^{(i,j)}$ recursively as
follows:
\[
\left\{
\begin{array}
[c]{lll}%
c^{(i,0)} & = & c^{(i)}\\
&  & \\
c^{(i,j+1)} & = & \sigma_{2^{i}}\left(  c^{(0,j)},\;\mathcal{S}c^{(1,j)}%
,\;\mathcal{S}^{2}c^{(2,j)},\;\ldots\;\mathcal{S}^{p}c^{(p,j)},\;\ldots
\right)
\end{array}
\right.
\]
Then $\lim_{j\rightarrow\infty}c^{(i,j)}$ exists and is given by
\[
\lim_{j\rightarrow\infty}c^{(i,j)}=\left\{
\begin{array}
[c]{ccc}%
0 & \text{if} & i>0\\
&  & \\
\sum_{p=0}^{\infty}\mathcal{S}^{p}c^{(p,q)} & \text{for all }q\geq0 & \text{if
}i=0
\end{array}
\right.
\]

\end{proposition}

\bigskip\ 

\begin{corollary}
Let $a$ and $b$ be two generic dyadic integers. Define the generic dyadic
integer
\[
c^{(i)}=b_{i}\diamond a
\]
Then, by using the construction given in the above proposition, the generic
product of $a$ and $b$ is given by
\[
a\cdot b=\lim_{j\rightarrow\infty}c^{(0,j)}%
\]

\end{corollary}

\bigskip

It is amusing and insightful to note that the above Boolean decomposition is
based on the repeated application of the following well known combinatorial
formula which acts as a bridge between characteristic $0$ and characteristic
$2$ algebraic structure, namely:

\bigskip

\begin{theorem}
Let $\mathbf{s}=s_{m-1}s_{m-2}\ldots s_{2}s_{1}s_{0}$ be a binary string of
length $m$, and let $\sigma_{k}\left(  y_{0},y_{1},y_{2},\ldots y_{m-1}%
\right)  $ be the $k$-th ($0\leq k<m$) elementary symmetric function in the
free Boolean ring $\mathbb{B}\left\langle y_{0},y_{1},y_{2},\ldots
y_{m-1}\right\rangle $. \ Then the Hamming weight $Wt\left(  \mathbf{s}%
\right)  $ of the string $\mathbf{s}$ is given by the following formula%
\[
Wt\left(  \mathbf{s}\right)  =\sum_{j=0}^{\left\lfloor \lg m\right\rfloor
-1}\sigma_{2^{j}}\left(  \mathbf{s}\right)  \cdot2^{j}\text{ ,}%
\]
where the symmetric function is first evaluated in the field of two elements
$\mathbb{F}_{2}$, and then interpreted as the integer $0$ or $1$ in
$\mathbb{Z}$.\cite{Knuth1}
\end{theorem}

\bigskip

Thus, the $\sigma_{2^{j}}$'s are the $j$-th order carries in the
multiplication algorithm. \ This is now made explicit in the following topdown
formulation of the multiplication algorithm.

\bigskip

\begin{proposition}
Let $Mat_{\infty}\left(  \mathbb{B}\left\langle \underline{\underline{x}%
}\right\rangle \right)  $ be the set of all Boolean matrices of the form%
\[
M=\left(  m_{ij}\right)  =\left(
\begin{array}
[c]{ccccc}%
\cdots & m_{03} & m_{02} & m_{01} & m_{00}\\
\ldots & m_{13} & m_{12} & m_{11} & 0\\
\ldots & m_{23} & m_{22} & 0 & 0\\
\cdots & m_{33} & 0 & 0 & 0\\
\vdots & \vdots & \vdots & \vdots & \vdots
\end{array}
\right)  \text{ ,}%
\]
where $m_{ij}\in\mathbb{B}\left\langle \underline{\underline{x}}\right\rangle
$ and $m_{ij}=0$ for $i>j$. \ Let $\Omega$ be the map defined by
\[%
\begin{array}
[c]{ccc}%
\Omega:Mat_{\infty}\left(  \mathbb{B}\left\langle \underline{\underline{x}%
}\right\rangle \right)   & \longrightarrow & Mat_{\infty}\left(
\mathbb{B}\left\langle \underline{\underline{x}}\right\rangle \right)  \\
M=\left(  \ldots,\gamma_{3},\gamma_{2},\gamma_{1},\gamma_{0}\right)   &
\mapsto & \overset{\mathstrut}{\left(
\begin{array}
[c]{ccccc}%
\cdots & \sigma_{2^{0}}\left(  \gamma_{3}\right)   & \sigma_{2^{0}}\left(
\gamma_{2}\right)   & \sigma_{2^{0}}\left(  \gamma_{1}\right)   &
\sigma_{2^{0}}\left(  \gamma_{0}\right)  \\
\ldots & \sigma_{2^{1}}\left(  \gamma_{2}\right)   & \sigma_{2^{1}}\left(
\gamma_{1}\right)   & \sigma_{2^{1}}\left(  \gamma_{0}\right)   & 0\\
\ldots & \sigma_{2^{2}}\left(  \gamma_{1}\right)   & \sigma_{2^{2}}\left(
\gamma_{0}\right)   & 0 & 0\\
\cdots & \sigma_{2^{3}}\left(  \gamma_{0}\right)   & 0 & 0 & 0\\
\vdots & \vdots & \vdots & \vdots & \vdots
\end{array}
\right)  }=\Omega\left(  M\right)
\end{array}
\]
where $\gamma_{j}$ denotes the $j$-th column of $M$, for $0\leq j<\infty$.
\end{proposition}

Let $a$ and $b$ be positive generic dyadic integers. \ Then%
\[
a\cdot b=\widetilde{P}\lim_{k\longrightarrow\infty}\Omega^{k}\left(
M_{0}\right)  \text{ ,}%
\]
where
\[
M_{0}=\left(
\begin{array}
[c]{ccccc}%
\cdots & a_{3}b_{0} & a_{2}b_{0} & a_{1}b_{0} & a_{0}b_{0}\\
\ldots & a_{2}b_{1} & a_{1}b_{1} & a_{0}b_{1} & 0\\
\ldots & a_{1}b_{2} & a_{0}b_{2} & 0 & 0\\
\cdots & a_{0}b_{3} & 0 & 0 & 0\\
\vdots & \vdots & \vdots & \vdots & \vdots
\end{array}
\right)  \text{ ,}%
\]
and where $P$ projects each matrix of $Mat_{\infty}\left(  \mathbb{B}%
\left\langle \underline{\underline{x}}\right\rangle \right)  $ onto its first row.

\bigskip

\begin{remark}
Please note that if the columns of the matrices of $Mat_{\infty}\left(
\mathbb{B}\left\langle \underline{\underline{x}}\right\rangle \right)  $ are
written in reverse order, then $Mat_{\infty}\left(  \mathbb{B}\left\langle
\underline{\underline{x}}\right\rangle \right)  $ becomes a ring of upper
triangular matrices.
\end{remark}

\bigskip

\begin{remark}
Please note that each entry $\sigma_{2^{i}}\left(  \gamma_{j}\right)  $ is
well defined because each column $\gamma_{j}$ is a positive generic integer.
\end{remark}

\bigskip

\begin{remark}
The sub-Boolean ring of $\mathbb{B}\left\langle \underline{\underline{x}%
}\right\rangle $ of all elementary symmetric Boolean functions is a free a
free Boolean ring with free basis $\sigma_{1},\sigma_{2},\sigma_{2^{2}}%
,\ldots,\sigma_{2^{\left\lfloor \lg n\right\rfloor }}$. \ Hence, as vector
spaces, $\dim_{\mathbb{F}_{2}}\left(  \mathbb{B}\left\langle \sigma_{2^{\ast}%
}\right\rangle \right)  =\left\lfloor \lg\dim_{\mathbb{F}_{2}}\left(
\mathbb{B}\left\langle \underline{\underline{x}}\right\rangle \right)
\right\rfloor $. \ Consequently, $\Omega^{k}\left(  M_{0}\right)  $ converges
exponentially fast to a matrix with only two non-zero rows.
\end{remark}

\end{document}